\documentstyle[12pt,amssymb,mysection]{article}
\voffset = - 2.3 truecm
\newcommand{\boxx}{\ \vrule width 2.0 mm height 2.0 mm  depth 0 mm }
\newcommand{\goth}{\frak}

\newcounter{theo}[section]
\newcounter{lemma}[section]
\newcounter{defi}[section]
\newcounter{rem}[section]
\newcounter{cor}[section]
\catcode`@=11
\catcode`@=12
\def\thetheo{\thesection.\arabic{theo}}
\def\thelemma{\thesection.\arabic{lemma}}
\def\thedefi{\thesection.\arabic{defi}}
\def\therem{\thesection.\arabic{rem}}
\def\thecor{\thesection.\arabic{cor}}
\newenvironment{theo}{%
\refstepcounter{theo} {THEOREM \thetheo.}\bgroup\rm}{\egroup}
\newenvironment{lemma}{%
\refstepcounter{lemma} {LEMMA \thelemma.}\bgroup\rm}{\egroup}
\newenvironment{defi}{%
\refstepcounter{defi} {DEFINITION \thedefi.}\bgroup\rm}{\egroup}
\newenvironment{rem}{%
\refstepcounter{rem} {REMARK \therem.} \bgroup\rm}{\egroup}
\newenvironment{cor}{%
\refstepcounter{cor} {COROLLARY \thecor.}\bgroup\rm}{\egroup}
\begin{document}
\ \vskip 3cm
\begin{minipage}{14.0cm}
{\bf  ON THE THEORY  OF MATRIX-VALUED FUNCTIONS\\
 BELONGING TO THE SMIRNOV CLASS\\[0.5cm]
V.E. Katsnelson and B. Kirstein
}
\end{minipage}
\vskip .5cm

%%%%%%%%%%%%%%%%%%%%%%%%%%%%%%%%%%%%% ABSTRACT %%%%%%%%%%%%%%%%%%%%%%%%%%%%%%
\begin{abstract}
A theory of matrix-valued functions from the matricial Smirnov class
${\goth N}_n^+({\Bbb D})$
is systematically developed. In particular, the maximum principle of
 V.I.Smirnov, inner-outer factorization, the Smirnov-Beurling characterization of outer functions and an analogue of Frostman's theorem are presented for
matrix-valued functions from the Smirnov class ${\goth N}_n^+({\Bbb D})$.
 We also consider a family $F_{\lambda} =F-\lambda I$  of functions
 belonging to the matricial Smirnov class which is indexed by a complex
 parameter $\lambda$.
We show that with the  exception of a ''very small'' set of such  $\lambda$
 the corresponding inner factor in the inner-outer factorization of
the function $F_{\lambda}$ is a Blaschke-Potapov product.

 The main goal of this paper is to provide users
 of analytic matrix-function theory with a standard source for references
related to  the matricial Smirnov class.
\end{abstract}
%%%%%%%%%%%%%%%%%%%%%%%%%%%%%%%%%%%%%%%%%%%%%%%%%%%%%%%%%%%%%%%%%%%%%%%%%%%%%%%
%%%%%%%%%%%%%%%%%%%%%%%% NOTATIONS %%%%%%%%%%%%%%%%%%%%%%%%%%%%%%%%%%%%%%%%%%%%
%%%%%%%%%%%%%%%%%%%%%%%%%%%%%%%%%%%%%%%%%%%%%%%%%%%%%%%%%%%%%%%%%%%%%%%%%%%%%%%
\vskip 0.4cm
\noindent {NOTATIONS : } \quad ${\Bbb C}$ \quad - \quad the complex plane.\\
\quad ${\Bbb T}:=\{ t \in {\Bbb C} \quad : \quad |t| = 1 \}$ \quad - \quad the unit circle.\\
\quad ${\Bbb D}:=\{ z \in {\Bbb C} \quad : \quad |z| < 1 \}$ \quad - \quad the unit disc.\\
\quad ${\goth B}_{{\Bbb T}}$ \quad - \quad the $\sigma$ - algebra of Borel
 subsets of ${\Bbb T}$.\\
\quad $m$ \quad - \quad  normalized Lebesgue measure on the measurable space
 $({\Bbb T}, {\goth B}_{{\Bbb T}})$.\\
\quad ${\Bbb C}^n$ \quad - \quad the $n$-dimensional complex space equipped
with the usual Euclidean norm, i.e., for $x=(\xi_1,\ldots,\xi_n)^{\top}$ we
 define
$ \| x \|_{{\Bbb C}^n} :=\Big\{\sum\limits_{k=1}^n |\xi_k|^2\Big\}^{1/2}.
$\\[0.1cm]
%\[ \| x \|_{{\Bbb C}^n} := \sqrt{\sum\limits_{k=1}^n |\xi_k|^2}. \]
${\goth M}_n$ \quad - \quad the set of all complex $n \times n$ matrices
equipped with the standard matrix norm, namely if $M \in {\goth M}_n$ then
$ \| M \| := \sup\limits_{x \in {\Bbb C}^n \setminus \{ 0 \}}
 \| Mx \|_{{\Bbb C}^n}\big / \| x \|_{{\Bbb C}^n}. $\\[0.2cm]
${\goth C}_n := \{ M \in {\goth M}_n \quad : \quad \| M \| \leq 1 \}$ \quad
- \quad the subset of all contractive matrices in ${\goth M}_n$.\\
\quad $I_n$ \quad - \quad the $n \times n$ unit matrix.\\
As usual for $r \in {\Bbb R}$ we set $r^+:=\max\{r,0\}$ and
$r^-:=\max\{-r,0\}$.\\
Hence, $r=r^+-r^-$ and $|r|=r^++r^-$. In particular, if $a \in (0,\infty)$,
then
\begin{eqnarray} \ln^+{a} = \max\{\ln{a},0\} \quad &, & \quad \ln^-{a} =
 \max\left\{\ln{\frac{1}{a}},0\right\}, \nonumber\\
\ln{a} = \ln^+{a}-\ln^-{a} \quad &, & \quad |\ln{a}| = \ln^+{a} + \ln^-{a}.
 \nonumber
\end{eqnarray}\\
If $A \in {\Bbb C}^{p \times q}$, then the symbol $A^{\top}$ stands for the
transposed matrix.\\[0.2cm]

%%%%%%%%%%%%%%%%%%%%%%%%%%%%%%%%%%%%%%%%%%%%%%%%%%%%%%%%%%%%%%%%%%%%%%%%%%%%%%%
%%%%%%%%%%%%%%%%%%%%%%%%%%%%%%% SECTION 0 %%%%%%%%%%%%%%%%%%%%%%%%%%%%%%%%%%%%%
%%%%%%%%%%%%%%%%%%%%%%%%%%%%%%%%%%%%%%%%%%%%%%%%%%%%%%%%%%%%%%%%%%%%%%%%%%%%%%

\setcounter{section}{-1}
\begin{minipage}{15.0cm}
\section{\hspace{-0.4cm}.\hspace{0.2cm} PREFACE}
\end{minipage}\\[-0.2cm]

 In this paper, we discuss various aspects of a class of matrix-valued
functions
which is named after V.I. Smirnov who introduced it for the scalar case in
his famous paper [Sm]. It should be mentioned that the scalar Smirnov class
also appeared in early papers of Doob (see e.g. [Doo1], [Doo2] and the
bibliographies in the monographs Collingwood and Lohwater [CoLo] and Noshiro
[No] which contain references to many other related works of Doob).
For a collection of basic facts on the Smirnov class
and the intimately related function spaces named after Nevanlinna and
Hardy we refer the reader to the monographs of P.L. Duren [Dur], J.B. Garnett
[G], K. Hoffman [Hoff], P. Koosis [Koo], I.I. Privalov [Pri] and
M. Rosenblum and J. Rovnyak [RoRo2]. These books concentrate more or less on
function-theoretic properties of functions belonging to some of the mentioned
classes. In the last two decades much progress has been made in clearing
up topological and functional-analytic questions connected with the structure
of the Smirnov class (see e.g. Yanagihara [Y1] - [Y10],
 Yanagihara and Kawase [YK],
Yanagihara and Nakamura [YN], Stoll [St1], [St2], Roberts [Rob],
 Roberts and Stoll [RoSt1], [RoSt2],
Mochizuki [Mo1], [Mo2], Helson [Hel2] - [Hel4], McCarthy [McC],
 Camera [Cam]).

A systematic study of the matricial Smirnov class was mainly promoted by the
 work of \mbox{D.Z. Arov.} In his paper [Ar1] on Darlington synthesis
 matricial generalization
of \mbox{V.I. Smirnov's} important maximum principle was used in an essential
way, namely with
its aid a powerful criterion for proving the $J$-contractivity of a meromorphic
matrix function was established. Moreover, D.Z. Arov's description of all
Darlington representations of a given (pseudocontinuable) Schur function
is based on the concept of denominators. {\sl A pair $[b_1,b_2]$ of inner
matrix-valued functions} of appropriate sizes is {\sl called a denominator of
 a given meromorphic matrix-valued function $f$} of bounded characteristic if
 $b_1 f b_2$ belongs to the matricial Smirnov class.

Nehari interpolation and generalized bitangential Schur - Nevanlinna - Pick
interpolation are other important problems which turned out to be closely
related with the matricial Smirnov class. This is an immediate consequence
of D.Z. Arov's work [Ar3] - [Ar9] (see also Nicolau [Nic1], [Nic2]). In his
investigations on the corresponding inverse problem D.Z. Arov introduced
particular subclasses of $J$-inner functions which are now called  the classes
of Arov-regular and Arov-singular $J$-inner functions. Here a $J$-inner
function $V$ is called Arov-singular if $V$ and $V^{-1}$ belong to the
matricial Smirnov class. Furthermore, a $J$-inner function $W$ is called
left Arov-regular (resp. right Arov-regular) if it does not contain
any nonconstant Arov-singular right (resp. left) divisors. D.Z. Arov (see
[Ar3] - [Ar7]) proved that each $J$-inner function $W$ admits (an essentially
unique) factorizations
\[ W = W_{l,r} \cdot W_{l,s} = W_{r,s} \cdot W_{r,r} \]
where the $J$-inner functions $W_{l,s}$ and $W_{r,s}$ are Arov-singular
whereas the $J$-inner functions $W_{l,r}$ and $W_{r,r}$ are left
Arov-regular and right Arov-regular, respectively. Furthermore, D.Z. Arov
proved that a $J$-inner function is a left (resp. right) resolvent
matrix of a completely indeterminate bitangential Schur - Nevanlinna - Pick
interpolation problem if and only if it is left Arov-regular (resp. right
Arov-regular). For several connections between left and right
Arov-regularity we refer the reader to the papers [Kats1], [Kats2] where
essential connections between left and right Blaschke - Potapov products
were established. In this way the first author (see [Kats3], [Kats4])
was led to a weighted approximation problems for pseudocontinuable functions
belonging to the Smirnov class. The papers [Kats1] -[Kats3] laid the basis for
the study of an inverse problem for Arov-singular $J$-inner functions which
was considered in [AFK7].
The papers [Ar2], [AFK1] - [AFK6] deal with several completion problems for
$J$-inner functions with particular emphasis on various subclasses of
$J$-inner functions (Smirnov type, inverse Smirnov type, Arov-singular type).
Using the concept of Arov-singularity and Arov-regularity of $J$-inner
functions and the approximation method created in [Kats3],  A. J. Kheifets
[Kh] answered a question of D. Sarason [Sar1] (see also [Sar2]) on exposed
points in the Hardy space $H^1({\Bbb D})$.
 Prediction theory for multivariate stationary sequences formed an important
source for the development of the theory of matrix-valued holomorphic
functions (see Wiener and Masani [WM1], [WM2], Helson and Lowdenslager [HL1],
[HL2], Rozanov [Roz1], [Roz2], Masani [Ma1] - [Ma4]). In particular, the
matricial Hardy class $H_n^2({\Bbb D})$ (see Definition \ref{du51} below)
 became
an essential tool. It turned out that the basic problems of prediction theory
could be reformulated as analytic problems for appropriate functions belonging
to the Hardy class $H_n^2({\Bbb D})$ . Using functional-analytic methods,
Beurling's inner-outer factorization was generalized to $H_n^2({\Bbb D})$
(see Masani [Ma2], Rozanov [Roz1]). Moreover classical results due to Szeg\"o
[Sz1] - [Sz3], Kolmogorov [Kol] and Krein [Kr] were extended to the
 multivariate case. Here, it  turned out (see Devinatz [De]) that the matrix
version of Szeg\"o's factorization theorem and other results due to Wiener and
 Masani [WM1], [WM2] and Helson and Lowdenslager [HL1], [HL2] are not so much
 generalizations
of Szeg\"o's classical results as consequences of it. An algebraic treatment
of this theory was given by Helson [Hel1].\\
Carrying on from the theory of matrix-valued functions belonging to the Hardy
class $H_n^2({\Bbb D})$,  we will study  various aspects of outer functions
 from the matricial Smirnov class in this paper. In particular,
 we will extend the theory of
 inner-outer factorization to the matricial Smirnov class. A central topic
in our investigations is to describe the situation where the inner factor
in the inner-outer factorization of a matrix-valued Smirnov class function
is a Blaschke - Potapov product.
Moreover, we will consider a family of functions belonging to the matricial Smirnov class
which is indexed by a complex parameter $\lambda$. Then it will be shown
that with exception of a "very small" set of such parameters $\lambda$
the corresponding inner factor in the inner-outer factorization of the
function $F_{\lambda}$ is a Blaschke - Potapov product. Our methods to
prove this use a matrix generalization of logarithmic potentials. In this way,
we obtain a generalization of a classical theorem of Frostman [Fr] (see
also Heins [Hei] and Rudin [Ru1], [Ru2]). It should be mentioned that it was
\mbox{Yu.~P. Ginzburg} who was a pioneer in matrix (and in operator) generalizations
of Frostman's results (see [Gi6] and [GiTa1] - [GiTa3]).\\[0.2cm]

%%%%%%%%%%%%%%%%%%%%%%%%%%%%%%%%%%%%%%%%%%%%%%%%%%%%%%%%%%%%%%%%%%%%%%%%%%%%%%%
%%%%%%%%%%%%%%%%%%%%%%%%%%%%%%% SECTION 1 %%%%%%%%%%%%%%%%%%%%%%%%%%%%%%%%%%%%%
%%%%%%%%%%%%%%%%%%%%%%%%%%%%%%%%%%%%%%%%%%%%%%%%%%%%%%%%%%%%%%%%%%%%%%%%%%%%%%%

\begin{minipage}{15.0cm}
\section {\hspace{-0.4cm}.\hspace{0.2cm} ON THE MATRICIAL NEVANLINNA AND
 SMIRNOV CLASSES }
\end{minipage}\\[-0.1cm]
\setcounter{equation}{0}

For $F : {\Bbb D} \to {\goth M}_n$ and  $r \in [0,1)$,  we define the
 function
$F_{[r]} : {\Bbb T} \to {\goth M}_n$ via $t \to F(rt)$.

\begin{defi}
\label{du11}
{\sl  A matrix-valued function $F: {\Bbb D} \to {\goth M}_n$ is said to
{\rm  belong to the matricial Nevanlinna class} ${\goth N}_n({\Bbb D})$
if $F$ is holomorphic in ${\Bbb D}$ and if the family
 $\left( \ln^+{\| F_{[r]} \|} \right)_{r \in [0,1)}$
is bounded in ${\cal L}^1(m)$, or more precisely, if}
\begin{equation} \sup\limits_{r \in [0,1)} \int\limits_{{\Bbb T}}
 \ln^+{\| F_{[r]}(t) \|} \ \ m(dt) < +\infty . \label{1.1}
\end{equation}
\end{defi}

\begin{rem}
\label{ru11a}
Let $F: {\Bbb D} \to {\goth M}_n$ be a matrix-valued function which is
holomorphic in ${\Bbb D}$. Then $F$ belongs to ${\goth N}_n({\Bbb D})$ if and
only if the (subharmonic) function $\ln^+{\| F \|}$ has a harmonic
majorant in ${\Bbb D}$.
\end{rem}

The definition of the Smirnov class ${\goth N}^+({\Bbb D})$ and of its
matricial analogue ${\goth N}_n^+({\Bbb D})$ are connected with the notion
 of uniform integrability. Since this notion is not used very often  we give
 the  definition.

{ DEFINITION OF UNIFORM INTEGRABILITY. }
{\sl Let $(\Omega,{\goth A},\mu)$ be a measure space.\\
Then the family $(f_{\alpha})_{\alpha \in A}$ belonging to ${\cal L}^1(\Omega,{\goth A},\mu;{\Bbb C})$ is called
{\rm uniformly integrable with respect to $\mu$}  if the following conditions
are satisfied} :\\
\[\hspace*{0.5cm}{\rm (i)}
 \qquad \qquad \qquad  \qquad \qquad \qquad
\sup\limits_{\alpha \in A} \int\limits_{\Omega} |f_{\alpha}(t)| \ \ \mu(dt) < + \infty .
\qquad \qquad \qquad \qquad \qquad \qquad \qquad
 \
 \]\\[0.2cm]
\hspace*{0.6cm}{\rm (ii) } \parbox[t]{14.0cm}{\sl For every
 $\epsilon \in (0,\infty)$ there exists a  $\delta \in (0,\infty)$
(which depends only on $\epsilon$) such that for all $\alpha \in A$ and for all
$\Delta \in {\goth A}$, with $\mu(\Delta)<\delta$,  the inequality
\[\textstyle
\int\limits_{\Delta} |f_{\alpha}(t)| \ \ \mu(dt) < \epsilon
  \mbox{\vspace*{-10.0cm}} \]
is fulfilled.} \\[0.2cm]

\begin{rem}
\label{ru11}
If $\mu(\Omega)<+\infty$ and if for each fixed $\delta \in (0,\infty)$
there exist an $N(\delta) \in {\Bbb N}$ and a sequence $(X_{k,\delta})_{k=1}^{N(\delta)}$
from $\Delta$ such that $\Omega=\bigcup\limits_{k=1}^{N(\delta)}\Omega_{k,\delta}$
and $\mu(\Omega_{k,\delta})\leq \delta$ for all $k \in \{1,2,\ldots,N(\delta)\}$,
 then a family of functions for which condition (ii) in the preceding
definition is fulfilled, automatically satisfies condition (i) .
Consequently, in the case of a finite measure space $(\Omega,{\goth A},\mu)$
condition (i) can be omitted in the definition of uniform integrability.
A special case of such a measure space is   the Lebesgue space
on ${\Bbb T}$,  where ${\goth A}$ is the $\sigma$ - algebra of Borel subsets
of ${\Bbb T}$ and $m$ is the normalized Lebesgue measure on ${\Bbb T}$.
\end{rem}
\ \\[0.3cm]
In the sequel we will repeatedly use the following theorem from measure theory
which goes back to G. Vitali [Vit] (see also [Ru3, p.133, Exercise 10]).

{ VITALI'S CONVERGENCE THEOREM.  } {\sl
 Let $(\Omega,{\goth A},\mu)$ be a finite
measure space (i.e., $\mu(\Omega)<\infty$). Let $(f_n)_{n \in {\Bbb N}}$ be a
sequence from ${\cal L}^1(\Omega,{\goth A},\mu;{\Bbb C})$ which is uniformly
integrable with respect to $\mu$ and converges $\mu$-a.e. to a Borel measurable
function $f: \Omega \to {\Bbb C}$.\\ Then $f \in {\cal L}^1(\Omega,{\goth A},\mu;{\Bbb C})$,
\[ \lim\limits_{n \to \infty} \int\limits_{\Omega} |f_n-f| \ d\mu = 0 \]
and
\[ \lim\limits_{n \to \infty} \int\limits_{\Omega} f_n \ d\mu = \int\limits_{\Omega} f \ d\mu . \]}

{ PROOF. } Let $\epsilon \in (0,\infty)$. In view of the uniform $\mu$-integrability
of $(f_n)_{n \in {\Bbb N}}$ there exists a number $\delta \in (0,\infty)$
such that for all $n \in {\Bbb N}$ and for all $\Delta \in {\goth A}$, which
satisfy $\mu(\Delta)<\delta$, the inequality
\begin{equation} \int\limits_{\Delta} |f_n| \ d\mu < \frac{\epsilon}{3} \label{1.2a}
\end{equation}
is satisfied. Since  $\mu(\Omega)<\infty$,  Egorov's Theorem guarantees
the existence of a set $B_{\delta} \in {\goth A}$ such that
\begin{equation} \mu(B_{\delta}) < \delta \label{1.3a}
\end{equation}
and
\begin{equation} \lim\limits_{n \to \infty} \sup\limits_{w \in \Omega \setminus B_{\delta}} |f_n(\omega) - f(\omega)| = 0 . \label{1.4a}
\end{equation}
Thus, there exists an $n_0 \in {\Bbb N}$ such that for all $n \geq n_0$ and
all $\omega \in \Omega \setminus B_{\delta}$ the inequality
\begin{equation} |f_n(\omega) - f(\omega)| < \frac{\epsilon}{3[1+\mu(\Omega)]} \label{1.5a}
\end{equation}
is satisfied. In view of $(\ref{1.2a})$ and $(\ref{1.3a})$ for $n \in {\Bbb N}$
we have
\begin{equation} \int\limits_{B_{\delta}} |f_n| \ d\mu < \frac{\epsilon}{3}.
 \label{1.6a}
\end{equation}
>From Fatou's Theorem and (\ref{1.6a}) we obtain
\begin{equation} \int\limits_{B_{\delta}} |f| \ d\mu \leq \mathop{\underline{\lim}}\limits_{n \to \infty} \int\limits_{B_{\delta}} |f_n| \ d\mu \leq \frac{\epsilon}{3}. \label{1.7a}
\end{equation}
Combining (\ref{1.5a}) - (\ref{1.7a}) we obtain  the estimate
\begin{eqnarray} \int\limits_{\Omega} |f_n-f| \ d\mu & = & \int\limits_{\Omega \setminus B_{\delta}} |f_n-f| \ d\mu + \int\limits_{B_{\delta}} |f_n-f| \ d\mu \nonumber\\
& \leq & \frac{\epsilon}{3[1+\mu(\Omega)]} \ \mu(\Omega \setminus B_{\delta}) + \int\limits_{B_{\delta}} |f| \ d\mu + \int\limits_{B_{\delta}} |f_n| \ d\mu \nonumber\\
& < & \frac{\epsilon}{3} + \frac{\epsilon}{3} + \frac{\epsilon}{3} = \epsilon \nonumber
\end{eqnarray}
for $n \geq n_0$. Thus,
\[ \lim\limits_{n \to \infty} \int\limits_{\Omega} |f_n-f| \ d\mu = 0 . \]
From this, all the remaining assertions follow immediately. \hfill $\boxx$

\begin{defi}
\label{du12a}
{\sl A function $\varphi: {\Bbb R} \to {\Bbb R}$ is called
{\rm strongly convex} if
 it has
the following properties}:\\[-0.1cm]
$$
\begin{array}{rl}
{\rm (\:i\:)}& \varphi\ {\sl is\ convex.}\\

{\rm(\hspace{0.04cm}ii\hspace{0.04cm})}& \varphi\
 {\sl is\ monotonically\ nondecreasing.}\\

{\rm(iii)}& \varphi\ {\sl takes\ its\ values\ in\ } [0,\infty).\\

{\rm (\hspace{0.01cm}iv\hspace{0.01cm})}&
 \displaystyle\lim\limits_{x \to \infty} \frac{\varphi(x)}{x} = \infty . \\

{\rm (\,v\,)} &
\parbox[t]{15cm}  {$ {\sl For\ some\ } c \in (0,\infty) {\sl \ there\
  exist\  constants\ } M
\in [0,\infty)
{\sl  \  and\  }
a \in {\Bbb R} {\sl \  such\  that\ }
\displaystyle \varphi(t+c) \leq M \cdot \varphi(t)
{\sl \ for \  all\ } t \in [a,\infty)$.}
\end{array}
$$
\end{defi}

If (v) holds for just one value of $c \in (0,\infty)$, say $c=c_0$, then by
 (ii) it holds for all $c \in (0,c_0)$. By iteration it holds for $c=nc_0$,
 $n \in {\Bbb N}$
and hence it holds for all $c \in (0,\infty)$.
%\ \\[0.5cm]

\begin{theo}
\label{tu11}
{ (de la Vall\'ee Poussin [LVP1], Nagumo [Na].)}{ \sl Let $(\Omega, {\goth A},
 \mu)$ be a (finite or infinite) measure space, and let
$(f_{\alpha})_{\alpha \in A}$
be a family of functions belonging to ${\cal L}^1(\Omega, {\goth A}, \mu; {\Bbb C})$.
In case $\mu(\Omega)=+\infty$, we assume also that
\[ \sup\limits_{\alpha \in A} \int\limits_{\Omega} |f_{\alpha}| \ d\mu < \infty . \]

{\rm (i)} \ \parbox[t]{13cm} {Suppose that there exists a function $\varphi: [0,\infty) \to [0,\infty)$
satisfying
\[ \lim\limits_{x \to +\infty} \frac{\varphi(x)}{x} = + \infty \]
and
\[ \sup\limits_{\alpha \in A} \int\limits_{\Omega} \varphi ( |f_{\alpha}| )\ d\mu < +\infty . \]
Then the family $(f_{\alpha})_{\alpha \in A}$ is uniformly integrable with respect
to $\mu$.\\[0.1cm]}

{\rm (ii)} \ \parbox[t]{13cm}{ Suppose that the family
 $(f_{\alpha})_{\alpha \in A}$ is uniformly
integrable with respect to $\mu$. Then there exists a strongly convex function
$\varphi: {\Bbb R} \to {\Bbb R}$ such that
\[ \sup\limits_{\alpha \in A} \int\limits_{\Omega} \varphi ( |f_{\alpha}|) \ d\mu < +\infty . \]}
}
\end{theo}

For a modern proof of Theorem \ref{tu11} we refer to [RoRo2, Theorem 3.10]
(see also Theorem 3.1.2 in [Ru2]). This modern proof based on  Vitali's
Convergence Theorem.

\begin{defi}
\label{du12}
{\sl A matrix-valued function $F: {\Bbb D} \to {\goth M}_n$ is said to
{\rm  belong to the matricial Smirnov class ${\goth N}_n^+({\Bbb D})$}
if $F$ is holomorphic in ${\Bbb D}$ and if the family $\left( \ln^+{\| F_{[r]}
 \|} \right)_{r \in [0,1)}$
is uniformly integrable with respect to the normalized Lebsgue measure $m$,
 i.e., if
for each $\epsilon \in (0,\infty)$ there exists a $\delta \in (0,\infty)$
(which depends only on $\epsilon$) such that for all $r \in [0,1)$ and for
all Borel subsets $\Delta$ of ${\Bbb T}$ satisfying $m(\Delta)<\delta$ the
inequality
\begin{equation} \int\limits_{\Delta} \ln^+{\| F_{[r]}(t) \|} \ \ m(dt) < \epsilon  \label{1.2}
\end{equation}
is fulfilled.}
\end{defi}

\begin{rem}
\label{ru12}
In view of Remark \ref{ru11}, each matrix-valued function $F \in {\goth N}_n^+({\Bbb D})$
satisfies condition (\ref{1.1}). Hence,
{\sl the matricial Smirnov class $ {\goth N}_n^+({\Bbb D})$ is a subclass
of the matricial Nevanlinna class ${\goth N}_n({\Bbb D})$}:
\begin{equation} {\goth N}_n^+({\Bbb D}) \subseteq {\goth N}_n({\Bbb D}) . \label{1.3}
\end{equation}
\end{rem}

For a matrix-valued function $F$ belonging to ${\goth N}_n({\Bbb D})$ we
denote by $\underline{F} : {\Bbb T} \to {\goth M}_n$ a boundary limit function
associated with $F$, i.e., $\underline{F}$ is a Borel measurable function and there
exists a Borel subset $\Delta_0$ of ${\Bbb T}$ satisfying $m(\Delta_0)=0$
such that for all $t \in {\Bbb T} \setminus \Delta_0$ we have
\[ \lim\limits_{r \to 1-0} F(rt) = \underline{F}(t) . \]
Observe that in view of Vitali's theorem a function $F \in {\goth N}_n^+({\Bbb D})$
satisfies
\begin{equation} \lim\limits_{r \to 1-0} \int\limits_{{\Bbb T}} \ln^+{\| F(rt) \|} \ \ m(dt) = \int\limits_{{\Bbb T}} \ln^+{\| \underline{F}(t) \|} \ \ m(dt) . \label{1.4}
\end{equation}
According to Fatou's theorem,
\begin{equation} \mathop{\underline{\lim}}\limits_{r \to 1-0} \int\limits_{{\Bbb T}} \ln^-{\| F(rt) \|} \ \ m(dt) \geq \int\limits_{{\Bbb T}} \ln^-{\| \underline{F}(t) \|} \ \ m(dt) . \label{1.5}
\end{equation}
(where equality does not hold in general). Hence,
\begin{equation} \mathop{\overline{\lim}}\limits_{r \to 1-0} \int\limits_{{\Bbb T}} \ln{\| F(rt) \|} \ \ m(dt) \leq \int\limits_{{\Bbb T}} \ln{\| \underline{F}(t) \|} \ \ m(dt) . \label{1.6}
\end{equation}

\begin{lemma}
\label{lu11}
{\sl A matrix-valued function $F : {\Bbb D} \to {\goth M}_n$ {\rm  belongs to
 the matricial
class ${\goth N}_n({\Bbb D})$ (resp. ${\goth N}_n^+({\Bbb D})$)} if and only
 if each
of its entries belongs to the scalar class ${\goth N}({\Bbb D})$
(resp. ${\goth N}^+({\Bbb D})$).}
\end{lemma}

As the determinant of a matrix is a polynomial of its elements and because each of the classes
${\goth N}({\Bbb D})$ and ${\goth N}^+({\Bbb D})$ is an algebra over ${\Bbb C}$ the following result
holds true.

\begin{lemma}
\label{lu12}
{\sl \parbox[t]{10cm}
{
{\rm (\hspace{1pt}i\hspace{1pt}) } If $F \in {\goth N}_n({\Bbb D})$,
 then $\mbox{det} F \in {\goth N}({\Bbb D})$.\\
{\rm (ii) } If $F \in {\goth N}_n^+({\Bbb D})$,
 then $\mbox{det} F \in {\goth N}^+({\Bbb D})$.}
}
\end{lemma}
\ \\

As a special case of (\ref{1.4}), (\ref{1.5}) and (\ref{1.6}) (corresponding
to the scalar case) we obtain for a function $F \in {\goth N}_n^+({\Bbb D})$
from part (ii) of Lemma \ref{lu12} that
\begin{equation} \lim\limits_{r \to 1-0} \int\limits_{{\Bbb T}}
 \ln^+{| \mbox{det} [F(rt)] |} \ \ m(dt) = \int\limits_{{\Bbb T}} \ln^+{| \mbox{det} [\underline{F}(t)] |} \ \ m(dt) , \label{1.7}
\end{equation}
\begin{equation} \mathop{\underline{\lim}}\limits_{r \to 1-0} \int\limits_{{\Bbb T}} \ln^-{| \mbox{det} [F(rt)] |} \ \ m(dt) \geq \int\limits_{{\Bbb T}} \ln^-{| \mbox{det} [\underline{F}(t)] |} \ \ m(dt) , \label{1.8}
\end{equation}
and, finally, that
\begin{equation} \mathop{\overline{\lim}}\limits_{r \to 1-0} \int\limits_{{\Bbb T}} \ln{| \mbox{det} [F(rt)] |} \ \ m(dt) \leq \int\limits_{{\Bbb T}} \ln{| \mbox{det} [\underline{F}(t)] |} \ \ m(dt) . \label{1.9}
\end{equation}

In the following we will use the {\sl Poisson kernel} $P: {\Bbb D} \times {\Bbb T} \to (0,\infty)$
which is defined by the formula
\[ P(z,t):=\frac{1-|z|^2}{|t-z|^2} . \]
 \ \\[-0.4cm]

\begin{theo}
\label{tu12}
{\sl
Let $F \in {\goth N}_n({\Bbb D})$ with $F \not\equiv 0$ and let  $u_F$ denote
the least harmonic majorant of $\log{\|F\|}$. Then the following statements are
 equivalent}:\\[-0.1cm]
$$
\begin{array}{rl}
{\rm(i)} &  F \in {\goth N}_n^+({\Bbb D}). \cr
\ & \ \cr
{\rm(ii)} &
 u_F(z) \leq \int\limits_{{\Bbb T}} \ln{\| \underline{F}(t) \|} \
\displaystyle  \frac{1-|z|^2}{|t-z|^2} \ m(dt) {\sl \
 \ for\  every\ }\displaystyle z \in {\Bbb D}.  \cr
\  & \ \cr {\rm (iii)} & \displaystyle \ln{ \| F(z) \|} \leq
\int\limits_{{\Bbb T}} \ln{\| \underline{F}(t) \|}
 \ \frac{1-|z|^2}{|t-z|^2} \ m(dt)
 {\sl \  \ for\  every\ }\displaystyle z \in {\Bbb D}. \cr
\  & \ \cr
{\rm(iv)} & \parbox[t]{14.0cm}{$ {\sl There \  exist \ a \ strongly \ convex \
 function } \  { \varphi:
 {\Bbb R} \to {\Bbb R}}
{\sl \  and \ a \ number\ } {r_0 \in (0,1)} {\sl \  such} \\[0.2cm]
 {\sl  that \  } \displaystyle
\   \sup\limits_{r \in [r_0,1)} \int\limits_{{\Bbb T}} \varphi \left( \ln{
 \| F_{[r]}(t) \|} \right) \ m(dt) < + \infty .$} \cr
\  & \ \cr
{\rm(v)}  &\parbox[t]{14.0cm}{$
  {\sl There \ exists \ a \ strongly \ convex \ function}\
\psi: {\Bbb R} \to {\Bbb R}
{\sl \  such \  that}\\[0.2cm]
\displaystyle \sup\limits_{r \in [0,1)} \int\limits_{{\Bbb T}} \psi \left( \ln^+{ \| F_{[r]}(t) \|} \right) \ m(dt) < + \infty . $}
\end{array}
$$
 \end{theo}
\ \\[-0.2cm]

{ PROOF. } Theorem \ref{tu12} can be proved by a slight modification of
the proof of Theorem 3.3.5 in [Ru2]. Here, Theorem \ref{tu11} plays an
essential role. \hfill $\boxx$

For further results on matrix-valued functions belonging to one of the classes
named after Nevanlinna, Smirnov and Hardy we refer the reader to
chapter 4 in [RoRo1].\\[0.2cm]

%%%%%%%%%%%%%%%%%%%%%%%%%%%%%%%%%%%%%%%%%%%%%%%%%%%%%%%%%%%%%%%%%%%%%%%%%%%%%%%
%%%%%%%%%%%%%%%%%%%%%%%%%%% SECTION 2 %%%%%%%%%%%%%%%%%%%%%%%%%%%%%%%%%%%%%%%%%
%%%%%%%%%%%%%%%%%%%%%%%%%%%%%%%%%%%%%%%%%%%%%%%%%%%%%%%%%%%%%%%%%%%%%%%%%%%%%%%

\begin{minipage}{15.5cm}
\section{\hspace{-0.4cm}. MATRIX FUNCTIONS OF THE SMIRNOV CLASS  AS
MUL\-TIPLES
 OF CONTRACTIVE MATRIX FUNCTIONS}
\end{minipage}\\[-0.1cm]
\setcounter{equation}{0}

\ \ \ Recall that a scalar function $e: {\Bbb D} \to {\Bbb C}$ is said to be
{\sl outer}
(in the sense of V.I. Smirnov) if there exist a unimodular constant $C \in {\Bbb T}$
and a function $w: {\Bbb T} \to [0,\infty)$ for which $\log{w}$ is
$m$-integrable such that for $z \in {\Bbb D}$ the relation
\begin{equation} e(z) = C \cdot \exp \left\{ \int\limits_{{\Bbb T}}
 \frac{t+z}{t-z} \ \ln{[w(t)]} \ \ m(dt) \right\} \label{2.1}
\end{equation}
holds true. Let ${\goth E}({\Bbb D})$ denote the class of all outer
 functions.
>From its definition it is obvious, that the class ${\goth E}({\Bbb D})$ is
 multiplicative:
If $e_1 , e_2 \in {\goth E}({\Bbb D})$, then $e_1 \cdot e_2 \in
 {\goth E}({\Bbb D})$.

The following statement is well-known (see e.g. Theorem 4.29 in [RoRo2]).

\begin{lemma}
\label{lu21}
{\sl Let $e: {\Bbb D} \to {\Bbb C}$ be some function.
Then the following statements are equivalent}:\\[-0.2cm]
$$
\begin{array}{rl}
{\rm (i) } &   e \in {\goth E}({\Bbb D}). \cr
\vspace*{-0.2cm} &\vspace*{-0.2cm} \cr
{\rm (ii) } &   e \in {\goth N}^+({\Bbb D}) \ , \ e \not\equiv 0$ \ and \
 $e^{-1} \in {\goth N}^+({\Bbb D}).\hspace{5.0cm}
\end{array}
$$
\end{lemma}
\\[-0.2cm]

In particular, a function $e$ of type (\ref{2.1}) belongs to the class
 ${\goth N}({\Bbb D})$
and, consequently, it possesses a boundary function $\underline{e}: {\Bbb T} \to {\Bbb C}$.
It is known that for almost all $t \in {\Bbb T}$ with respect to $m$,
\begin{equation} |\underline{e}(t)|=w(t) . \label{2.2}
\end{equation}
As the function $\ln{|e|}$ is harmonic in ${\Bbb D}$
we obtain
\[ \int\limits_{{\Bbb T}} \ln{|e(rt)|} \ \ m(dt) = \ln{|e(0)|} = \int\limits_{{\Bbb T}} \ln{[w(t)]} \ \ m(dt). \]
for $r \in [0,1)$.
Consequently, if $e \in {\goth E}({\Bbb D})$, then for $r \in [0,1)$ we obtain
\begin{equation} \int\limits_{{\Bbb T}} \ln{|e(rt)|} \ \ m(dt) = \int\limits_{{\Bbb T}} \ln{|\underline{e}(t)|} \ \ m(dt) . \label{2.3}
\end{equation}

Let us recall the following useful characterization of outer functions
 (see e.g.
Corollaries 4.16 and 4.17 in [RoRo2]).

\begin{lemma}
\label{lu21a}
{\sl
Let $e \in {\goth N}({\Bbb D})$  but $e \not\equiv 0$. Then the following
statements are\\ equivalent}:
$$
\begin{array}{rl}
 {\rm (i) } &  e  {\sl \ is \ outer}.\cr

{\rm (ii) } & {\sl  For \  all\ } z \in {\Bbb D}, \quad \displaystyle
 \ln{|e(z)|} = \int\limits_{{\Bbb T}} \mbox{Re} \frac{t+z}{t-z}
 \ \ln{| \underline{e}(t)|} \ m(dt) . \cr

{\rm (iii) } & {\sl  There \ is \ a \ } z_0 \in   {\Bbb D} {\sl \  such \
 that\ }
\displaystyle \ln{|e(z_0)|} = \int\limits_{{\Bbb T}}\mbox{Re}
 \frac{t+z_0}{t-z_0} \
 \ln{| \underline{e}(t)|} \ m(dt) .
\cr

{\rm (iv) }  &  {\sl If\ } {h \in {\cal N}^+({\Bbb D})}
 {\sl \  satisfies\ }\
{ |\underline{h}(t)| \leq |\underline{e}(t)|}\
{\sl \ for\  almost \ all \ } t {\sl \in } {\Bbb T}
 \ {\sl with\  respect \ to  }\   m ,{\sl \ then}\cr
\vspace*{-0.2cm} & \vspace*{-0.2cm}  \cr
 \ &{ for \ all \ } {z \in {\Bbb D}}\  {\sl  the\ inequality\ }
\  |h(z)| \leq |e(z)| {\sl \ holds \ true.}\cr
\vspace*{-0.0cm} &\vspace*{-0.0cm} \cr
{\rm (v) } &  {\sl  If\ } z_0 \in {\Bbb D}
{\sl \ and \ if\ } h \in {\cal N}^+({\Bbb D})
{\sl \ satisfies\ } |\underline{h}(t)| \leq |\underline{e}(t)|
{\sl \ for\  {\sl almost} \ all \ } t \in {\Bbb T}
\ {\sl with\  respect \ to  }\    \cr
\ & m , {\sl  then  \ the  \
 inequality\ } |h(z_0)| \leq |e(z_0)|{\sl \ holds \ true.}
\end{array}
$$
\end{lemma}
\ \\

Observe that conditions (iii) and (v) of Lemma \ref{lu21a} are usually
used with the choice $z_0=0$.

In the proof of Lemma \ref{lu23} and also in further considerations  we will use the
following result which goes back to V.I. Smirnov [Sm].
\ \\[-0.1cm]

{\rm THE MAXIMUM PRINCIPLE OF V.I. SMIRNOV. } {\sl
 Let $f \in {\goth N}^+({\Bbb D})$
be such that its boundary function $\underline{f}$ is $m$-essentially
 bounded.
Then $f$ is bounded in the unit disc and satisfies
\[ \sup\limits_{z \in {\Bbb D}} |f(z)| =
{\mbox{\rm ess}} \sup\limits_{\mbox{\hspace{-0.6cm}}t\in{\Bbb T}}
 | {\underline{f}}(t)| \]}

This result can be generalized to the matrix case.\\[0.1cm]

{\rm THE MAXIMUM PRINCIPLE OF V.I. SMIRNOV FOR MATRIX FUNCTIONS. }\\
{\sl Let $F \in {\goth N}_n^+({\Bbb D})$ be such that its boundary function
 $\underline{F}$ satisfies
${\mbox{\rm ess}} \sup\limits_{\mbox{\hspace{-0.6cm}}t\in{\Bbb T}}
 \| \underline{F}(t) \| < \infty . $\\
Then $F$ is bounded in the unit disc and satisfies
\[ \sup\limits_{z \in {\Bbb D}} \| F(z) \| =
{\mbox{\rm ess}} \sup\limits_{\mbox{\hspace{-0.6cm}}t\in{\Bbb T}}
  \| \underline{F}(t) \| . \]
}
\ \\[-0.8cm]

{\rm  PROOF. } Let $F=(F_{j,k})_{j,k=1}^n$, and fix the indices
 $j,k \in \{ 1,\ldots,n \}$.
In view of the inequality
 $|F_{j,k}(z)| \leq \| F(z) \|   \  (z \in {\Bbb D}) \ \ $, then
$F_{j,k} \in {\goth N}^+({\Bbb D})$ and
\[{\mbox{\rm ess}} \sup\limits_{\mbox{\hspace{-0.6cm}}t\in{\Bbb T}}
 | \underline{F_{j,k}}(t)| \leq
{\mbox{\rm ess}} \sup\limits_{\mbox{\hspace{-0.6cm}}t\in{\Bbb T}}
 \| \underline{F}(t) \| . \]
According to the maximum principle for scalar functions we then have
\[ \sup\limits_{z \in {\Bbb D}} |F_{j,k}(z)| < +\infty . \]
Hence,
\[ \sup\limits_{z \in {\Bbb D}} \| F(z) \| < +\infty . \]
The bounded holomorphic matrix-valued function $F$ admits the Poisson integral
representation
\[ F(z) = \int\limits_{{\Bbb T}} \underline{F}(t) \cdot P(z,t) \ m(dt) \quad , \quad z \in {\Bbb D} . \]
Therefore, by the integral version of the triangle inequality,  we obtain
\[ \| F(z) \| \leq \int\limits_{{\Bbb T}} \| \underline{F}(t) \| \cdot P(z,t) \ m(dt) \quad , \quad z \in {\Bbb D} . \]
But this in turn implies the inequality
$ \| F(z) \| \leq
{\mbox{\rm ess}} \sup\limits_{\mbox{\hspace{-0.6cm}}t\in{\Bbb T}}
 \| \underline{F}(t) \| \ \ (z \in {\Bbb D}).$ \hfill \boxx

Since  the function $\ln^+{\|F\|}$ is subharmonic  for an analytic
 matrix-valued function $F$ the following result is true.

\begin{lemma}
\label{lu22}
{\sl
Let $F \in {\goth N}_n^+({\Bbb D})$. Then for all $z \in {\Bbb D}$ the inequality
\begin{equation} \| F(z) \| \leq \exp \left\{ \int\limits_{{\Bbb T}} P(z,t) \ln{\|\underline{F}(t)\|} \ \ m(dt) \right\} \label{2.4}
\end{equation}
holds true.}
\end{lemma}

For a proof of Lemma \ref{lu22} we refer to Theorem 3.13 in [RoRo2].

Clearly, the maximum principle of V.I. Smirnov is a consequence of
inequality (\ref{2.4}).

\begin{defi}
\label{du21}
{\sl
The set ${\goth S}_{n \times n}({\Bbb D})$ of all holomorphic matrix-valued
functions $S: {\Bbb D} \to {\goth C}_n$ is called }{\rm
 the $n \times n$ Schur class.}
\end{defi}

\begin{lemma}
\label{lu23}
{\sl
A matrix-valued function $F: {\Bbb D} \to {\goth M}_n$ belongs to the Smirnov
class ${\goth N}_n^+({\Bbb D})$ if and only if it admits a representation
of the form
\begin{equation} F=\frac{1}{d} \cdot \Phi ,\label{2.5}
\end{equation}
where $\Phi \in {\goth S}_{n \times n}({\Bbb D})$ and $d$ is an  outer
function which belongs to ${\goth S}({\Bbb D}) $.}
\end{lemma}

PROOF. \  {\rm I. } Suppose that $F$ admits a representation of the form
 (\ref{2.5}).
Then $\Phi \in {\goth N}_n^+({\Bbb D})$ and,
as $d$ is outer, we have $d^{-1} \in {\goth N}^+({\Bbb D})$. Thus, as
 ${\goth N}^+({\Bbb D})$
is an algebra over ${\Bbb C}$, we get
 $\Phi \cdot d^{-1} \in {\goth N}_n^+({\Bbb D})$.\\[0.1cm]
{\rm II. } Suppose that $F \in {\goth N}_n^+({\Bbb D})$. We can assume that
 $F$ is not
the null function in ${\Bbb D}$. Then $\ln{\|\underline{F}\|}$ is $m$-integrable.
We define $d: {\Bbb D} \to {\Bbb C}$ via
\[ d(z) := \exp \left\{ - \int\limits_{{\Bbb T}} \frac{t+z}{t-z} \ln{\|\underline{F}\|} \ \ m(dt) \right\} . \]
Then, from our earlier considerations
 (see (\ref{2.1}) - (\ref{2.4})), it is clear
that $d$ is a scalar outer function and that the corresponding boundary
 function
$\underline{d}$ satisfies
\begin{equation} |\underline{d}(t)| = \|\underline{F}(t)\|^{-1} \label{2.6}
\end{equation}
for almost all $t \in {\Bbb T}$ with respect to $m$.
 Now define $\Phi : {\Bbb D} \to {\goth M}_n$ via
\begin{equation} \Phi(z) := d(z) \cdot F(z) . \label{2.7}
\end{equation}
Then, since $F \in {\goth N}_n^+({\Bbb D}), d \in {\goth N}^+({\Bbb D})$
 and  ${\goth N}^+({\Bbb D})$ is an algebra over ${\Bbb C}$, we see that
\begin{equation} \Phi \in {\goth N}_n^+({\Bbb D}) . \label{2.8}
\end{equation}
>From (\ref{2.6}) and (\ref{2.7}) we get
\begin{equation} \| \underline{\Phi}(t) \| =1 \label{2.9}
\end{equation}
for almost all $t \in {\Bbb T}$ with respect to $m$.
 Finally, in view of (\ref{2.8}) and (\ref{2.9}),
the maximum prin\-ciple of V.I. Smirnov implies that for $z \in {\Bbb D}$
we obtain $\| \Phi(z) \| \leq 1$.
Thus, $\Phi \in {\goth S}_{n \times n}({\Bbb D})$. \hfill $\boxx$\\[0.2cm]

%%%%%%%%%%%%%%%%%%%%%%%%%%%%%%%%%%%%%%%%%%%%%%%%%%%%%%%%%%%%%%%%%%%%%%%%%%%%%%%
%%%%%%%%%%%%%%%%%%%%%%%%%%%%%% SECTION 3 %%%%%%%%%%%%%%%%%%%%%%%%%%%%%%%%%%%%%%
%%%%%%%%%%%%%%%%%%%%%%%%%%%%%%%%%%%%%%%%%%%%%%%%%%%%%%%%%%%%%%%%%%%%%%%%%%%%%%%

\begin{minipage}{14.0cm}
\section{\hspace{-0.4cm}. OUTER MATRIX-VALUED FUNCTIONS}
\end{minipage}\\[-0.1cm]
\setcounter{equation}{0}

The main goal of this section is to discuss outer matrix-valued
functions which belong to the Smirnov class ${\goth N}_n^+({\Bbb D})$.
The needs of prediction theory of multivariate stationary stochastic processes
initiated an intensive study of matrix-valued outer functions belonging to the
Hardy class $H_n^2({\Bbb D})$ (see Definition \ref{du51} below) which is a
 subclass of ${\goth N}_n^+({\Bbb D})$.
The formula for the best predictor of a multivariate stationary stochastic
process of a given time in terms of its past depends in an essential manner
on a particular outer matrix-valued function belonging to $H_n^2({\Bbb D})$
(see Wiener and Masani [WM1], [WM2], Helson and Lowdenslager [HL1], [HL2],
Rozanov [Roz1], [Roz2], Masani [Ma1] - [Ma4] and for operator-valued
 generalizations also
Devinatz [De], Helson [Hel1], Sz.-Nagy and Foias [SZNF], Nikolskii [Nik2]).

\begin{defi}
\label{du31}
{\sl A matrix-valued function $E: {\Bbb D} \to {\goth M}_n$ is called
{\rm  outer}
(in the sense of V.I. Smirnov) if $E \in {\goth N}_n^+({\Bbb D})$ and
$\det{E}$ is outer. The class of all $n \times n$ matrix-valued outer
functions will be denoted by ${\goth E}_n({\Bbb D})$.
}
\end{defi}

If $E \in {\goth E}_n({\Bbb D})$ then, in particular, for all $z \in {\Bbb D}$
we have
\[ \det{[E(z)]} \not= 0 . \]

Definition \ref{du31} is clearly an immediate generalization of the notion of
a scalar outer function. This definition of an outer matrix-valued function
enables us to avoid the study of the question of a matricial analogue of
formula (\ref{2.1}).

\begin{rem}
The class ${\goth E}_n({\Bbb D})$ is multiplicative:
If $E_1,E_2 \in {\goth E}_n({\Bbb D})$ then $E_1 \cdot E_2 \in {\goth E}_n({\Bbb })$.
\end{rem}

\begin{rem}
\label{ru31a}
Let $E \in {\goth E}_n({\Bbb D})$. Then $E^{\top} \in {\goth E}_n({\Bbb D})$.
\end{rem}
\ \\[-0.2cm]

\begin{theo}
\label{tu31}
{(Determinant characterization of outer matrix-valued
 functions) }

 {\rm  (i) } {\sl  Let $E \in {\goth E}_n({\Bbb D})$. Then \ $\det{[E(z)]}
 \not=0$ \ for
all $z \in {\Bbb D}$ and \ $E^{-1} \in {\goth N}_n^+({\Bbb D})$.

{\rm (ii) } \parbox[t]{13.0cm}{Let $E$ be a function from ${\goth N}_n^+({\Bbb D})$ for which $\det{E}$
never vanishes in ${\Bbb D}$ and $E^{-1}$ belongs to ${\goth N}_n^+({\Bbb D})$.
Then $E \in {\goth E}_n({\Bbb D})$.}}
\end{theo}

 PROOF.
{\rm (i) } According to the rule for computing the inverse matrix
we have the representation
\begin{equation} E^{-1} = \frac{1}{\det{E}} \cdot A \label{3.1}
\end{equation}
where $A: {\Bbb D} \to {\goth M}_n$ is a matrix-valued function the entries of
 which
are polynomials of the elements of matrix $E$ (namely, the cofactors of the
corresponding elements). Since the class ${\goth N}^+({\Bbb D})$ is an algebra
 over ${\Bbb C}$, each entry of $A$ belongs to ${\goth N}^+({\Bbb D})$.
 Hence, $A \in {\goth N}_n^+({\Bbb D})$. From the fact that
 $E \in {\goth E}_n({\Bbb D})$ and Lemma \ref{lu21}  it then  follows that
 $(\det{E})^{-1} \in {\goth N}^+({\Bbb D})$, and thus in view of (\ref{3.1}),
 $E^{-1} \in {\goth N}_n^+({\Bbb D})$.
Hence, (i) is proved. \\[0.1cm]
{\rm (ii) } By Lemma \ref{lu12},  $\det{E} \in {\goth N}^+({\Bbb D})$ and
$\det{(E^{-1})} \in {\goth N}^+({\Bbb D})$. Therefore, the function
 $\det{E}$ satisfies condition (ii) in Lemma \ref{lu21}.
Thus,  $\det{E} \in {\goth E}({\Bbb D})$, and so, in view of Definition
 \ref{lu31},
 $E \in {\goth E}_n({\Bbb D})$.
Hence (ii) is proved. \hfill $\boxx$

The following result supplements the statement of Lemma \ref{lu23}.

\begin{lemma}
{\sl
\label{lu31}
Let $E \in {\goth E}_n({\Bbb D})$. Then $E$ has a representation of the form
\begin{equation} E=\frac{1}{d} \cdot C ,\label{3.2}
\end{equation}
where $C \in {\goth S}_{n \times n}({\Bbb D}) \cap {\goth E}_n({\Bbb D})$ and
$d \in {\goth E}({\Bbb D})$.}
\end{lemma}

 PROOF. \  In view of Lemma \ref{lu23}, the function $E$ has a representation
of the form
\[ E=\frac{1}{d} \cdot C ,\]
where $C \in {\goth S}_{n \times n}({\Bbb D})$ and $d \in {\goth E}({\Bbb D})$.
Lemma \ref{lu21} guarantees that
 $d^{-1} \in {\goth N}^+({\Bbb D})$. Since $E \in {\goth E}_n({\Bbb D})$,
it follows from Theorem \ref{tu31} that $E^{-1} \in {\goth N}_n^+({\Bbb D})$.
Therefore, as ${\goth N}^+({\Bbb D})$ is an algebra over ${\Bbb C}$,
 from \ $C^{-1} = d^{-1} \cdot E^{-1}$
we see that $C^{-1} \in {\goth N}_n^+({\Bbb D})$. Thus, as
$C \in {\goth S}_{n \times n}({\Bbb D}) \subseteq {\goth N}_n^+({\Bbb D})$
it follows from
Theorem \ref{tu31} that $C \in {\goth E}_n({\Bbb D})$. \qquad \hfill $\boxx$

Let us recall the following notion.

\begin{defi}
\label{du32}
{\sl The class $H_n^{\infty}({\Bbb D})$}  consists of all matrix-valued
functions  $F: {\Bbb D} \to {\goth M}_n$ which are holomorphic and bounded in
 ${\Bbb D}$,
i.e.,
\begin{equation} \sup\limits_{z \in {\Bbb D}} \| F(z) \| < \infty . \label{3.3}
\end{equation}
\end{defi}

\begin{theo}
\label{tu32}
{\rm (i) } {\sl Let $E \in {\goth E}_n({\Bbb D})$. Then there exists a sequence
$(F_k)_{k \in {\Bbb N}}$ from $H_n^{\infty}({\Bbb D})$ with the following
 properties:\\[0.2cm]
%%%%
\hspace*{1.0cm}$(\alpha)$ \parbox[t]{14cm}{ For almost
 all $t \in {\Bbb T}$ with respect to $m$, \quad
$ \displaystyle
 \lim\limits_{k \to \infty} \underline{E}(t) \cdot \underline{F_k}(t) = I_n .
$}\\[0.2cm]
\hspace*{1.0cm}$(\beta)$ \parbox[t]{14cm}{ The family $\left(\ln^+{\|\underline{F_k}\|}\right)_{k \in {\Bbb N}}$ is
uniformly integrable with respect to {\rm m}.}\\[0.2cm]
\hspace*{1.0cm}$(\gamma)$ \parbox[t]{14cm}{There exists a Borel subset $B_0$
 of ${\Bbb T}$ with $m(B_0)=0$
such that for all $k \in {\Bbb N}$ and all $t \in {\Bbb T} \setminus B_0$
 the inequality
$\displaystyle
 \| \underline{E}(t) \cdot \underline{F_k}(t) \| \leq 1
$
holds true.}\\[0.3cm]
\hspace*{0.6cm}{\rm (ii) } Let $E \in {\goth N}_n^+({\Bbb D})$ be such that there exists a sequence
$(F_k)_{k \in {\Bbb N}}$ belonging to $H_n^{\infty}({\Bbb D})$\\
\hspace*{1.2cm} satisfying
the above conditions $(\alpha)$ and $(\beta)$. Then $E \in {\goth E}_n({\Bbb D})$.}
\end{theo}
\ \\[0.3cm]
\begin{rem}
\label{ru31}
Theorem \ref{tu32} expresses in some sense a Smirnov class genera\-lization of that
characterization of the property that a function is outer which is formulated
in terms of the shift-invariant subspace generated by this function. Sometimes
the approximation property contained in Theorem \ref{tu32} is called weak
invertibility of the function $E$ (see [Sh] or [Nik1,Ch.2]). For the spaces
$H_n^{\infty}({\Bbb D})$ or $H_n^2({\Bbb D})$ this approximation property
(weak invertibility) will be often used for defining the notion "outer function".
Observe that in the scalar case $(n=1)$ it was already shown by V.I. Smirnov [Sm]
that for an outer function $e$ the linear subspace $e \cdot H^2({\Bbb D})$ is
dense in $H^2({\Bbb D})$. Concerning several generalizations of this result of
V.I. Smirnov we refer the reader to chapter 2 in [Nik1] (in particular, see
Theorem 3 in Section 2.2. ).
\end{rem}

 PROOF OF THEOREM \ref{tu32}.  {\rm (i) } Since $E$ is a matrix-valued outer
function, Theorem \ref{tu31} guarantees that
 $E^{-1} \in {\goth N}_n^+({\Bbb D})$.
We fix a boundary function $\underline{E}$ of $E$ such $\det{[\underline{E}(t)]} \not= 0$
for $t \in {\Bbb T}$. Then for $k \in {\Bbb N}$ we define $w_k : {\Bbb T} \to (0,\infty)$ via
\begin{equation} w_k(t) := \left\{ \begin{array}{ll}
\ \ \ \ \ \, 1 \ &, \mbox{ if} \quad \| \underline{E}^{-1}(t) \| < k \\
\ \ \\
\displaystyle \frac{1}
{\displaystyle\| \underline{E}^{-1}(t) \|} \
 &, \mbox{ if} \quad \| \underline{E}^{-1}(t) \| \geq k  .
\end{array} \right. \label{3.4} \end{equation}
Clearly
\begin{equation}
 0 < w_1(t) \leq w_2(t) \leq w_3(t) \leq \ldots \leq 1 \label{3.5}
\end{equation}
for $t \in {\Bbb T}$ and
\begin{equation} \lim\limits_{k \to \infty} w_k(t) = 1 . \label{3.6}
\end{equation}
From (\ref{3.5}) we see that  the inequality
\begin{equation} w_1(t) \geq \| \underline{E}^{-1}(t) \|^{-1} \label{3.7}
\end{equation}
holds for $t \in {\Bbb T}$. Since
 $E^{-1} \in {\goth N}_n^+({\Bbb D})$, we infer that
\begin{equation} \ln{[\| \underline{E}^{-1} \|^{-1}]} \in {\cal L}^1({\Bbb T}, {\goth B}_{{\Bbb T}}, m; {\Bbb C}) . \label{3.8}
\end{equation}
From (\ref{3.5}) - (\ref{3.7}) we obtain
\begin{equation} \int\limits_{{\Bbb T}} \ln{[w_k(t)]} \ m(dt) > -\infty . \label{3.9}
\end{equation}
Hence, for $k \in {\Bbb N}$ the function $\varphi_k : {\Bbb D} \to {\Bbb C}$
which is given by
\[ \varphi_k(z) := \exp\left\{ \int\limits_{{\Bbb T}} \frac{t+z}{t-z} \ln{[w_k(t)]} \ m(dt) \right\} \]
is well-defined. Moreover from its definition it is clear that $\varphi_k \in {\goth N}^+({\Bbb D})$
(or more precisely, that $\varphi_k$ is even outer). In view of (\ref{3.5}) and
(\ref{3.6}) the monotone convergence theorem guarantees that
\begin{equation} \lim\limits_{k \to \infty} \varphi_k(z) = 1 \ , \qquad z \in {\Bbb D} . \label{3.10}
\end{equation}
Since $|\underline{\varphi_k}(t)|=w_k(t)$ for almost all $t \in {\Bbb T}$
with respect to $m$,
formula (\ref{3.6}) yields
\[ \lim\limits_{k \to \infty} |\underline{\varphi_k}(t)| = 1 \]
for almost all $t \in {\Bbb T}$  with respect to $m$.
 In view of (\ref{3.5}) and
(\ref{3.6}),  another application of the monotone convergence theorem gives us
\begin{equation} \lim\limits_{k \to \infty} \int\limits_{{\Bbb T}} |\underline{\varphi_k}(t)|^2 \ m(dt) = \lim\limits_{k \to \infty} \int\limits_{{\Bbb T}} [w_k(t)]^2 \ m(dt)
= \int\limits_{{\Bbb T}} 1 \ dm = 1 . \label{3.11}
\end{equation}
For $k \in {\Bbb N}$, we have
\begin{equation} \int\limits_{{\Bbb T}} |\underline{\varphi_k}(t)-1|^2 \ m(dt) = \int\limits_{{\Bbb T}} |\underline{\varphi_k}(t)|^2 \ m(dt) - 2 \Re{[\varphi_k(0)]} + 1 . \label{3.12}
\end{equation}
Combining (\ref{3.10}) - (\ref{3.12}) it follows that
\begin{equation} \lim\limits_{k \to \infty} \int\limits_{{\Bbb T}} |\underline{\varphi_k}(t)-1|^2 \ m(dt) = 0 . \label{3.13}
\end{equation}
In view of (\ref{3.13}), the F. Riesz - Fischer theorem yields a subsequence $(\underline{\varphi_{l_k}})_{k \in {\Bbb N}}$
of $(\underline{\varphi_k})_{k \in {\Bbb N}}$ such that
\begin{equation} \lim\limits_{k \to \infty} \underline{\varphi_{l_k}} (t) = 1 \label{3.14}
\end{equation}
for almost all $t \in {\Bbb T}$ with respect to $m$. Let $k \in {\Bbb N}$ and
  set
\begin{equation} F_k:=E^{-1} \cdot \varphi_{l_k} . \label{3.15}
\end{equation}
Then, since $E^{-1} \in {\goth N}_n^+({\Bbb D})$ and
$ \varphi_{l_k} \in {\goth N}^+({\Bbb D})$,
we get $F_k \in {\goth N}_n^+({\Bbb D})$. Thus as $|\underline{\varphi_{l_k}}|
=w_{l_k}$
almost everywhere with respect to $m$ it follows from (\ref{3.15}) and
 (\ref{3.4}) that
\[ \| \underline{F_k} (t) \| = w_{l_k}(t) \cdot \| \underline{E}^{-1}(t) \|
 \leq l_k \]
for almost all $t \in {\Bbb T}$ with respect to $m$.
Thus, the maximum principle of V.I. Smirnov implies that
 $\| F_k(z) \| \leq l_k$ for all $z \in {\Bbb D}$ . Consequently,
 $F_k \in H_n^{\infty}({\Bbb D})$.
From (\ref{3.15}) it follows that
\begin{equation} E \cdot F_k = \varphi_{l_k} \cdot I_n . \label{3.16}
\end{equation}
From (\ref{3.5}) we obtain
\[ |\underline{\varphi_{l_k}}(t)|=w_{l_k}(t) \leq 1  \]
and hence since $\varphi_{l_k} \in {\goth N}^+({\Bbb D})$,
 the maximum principle of
V.I. Smirnov guarantees that that
\begin{equation} |\varphi_{l_k}(z)| \leq 1 , \qquad z \in {\Bbb D} . \label{3.17}
\end{equation}
Thus, combining (\ref{3.16}) and (\ref{3.17}) we see that $(\gamma)$ is fulfilled.\\
Moreover, from (\ref{3.16}) and (\ref{3.14}) we get that $(\alpha)$ is
 satisfied.\\
For almost all $t \in {\Bbb T}$ with respect to $m$
we have $|\underline{\varphi_{l_k}}(t)| \leq 1$
and, consequently, in view of (\ref{3.15}),  the inequality
\[ \ln^+{\|\underline{F_k}(t)\|} \leq \ln^+{\|\underline{E}^{-1}(t)\|} \]
holds  for almost all $t \in {\Bbb T}$ with respect to $m$. Hence,
 the family $(\ln^+{\|\underline{F_k}(t)\|})_{k \in {\Bbb N}}$
has an $m$-integrable majorant. This implies that $(\beta)$ is fulfilled.\\
Part (i) of Theorem \ref{tu32} is now proved.

Before proving part (ii) of Theorem \ref{tu32} we recall the following result
(see [WM1, Lemma 3.12]).\\[-0.1cm]

THE GENERALIZED MINKOWSKI INEQUALITY. {\sl  Let $(\Omega,{\goth A},P)$ be a
probability space and let $M: \Omega \to {\goth M}_n$ be a $P$-integrable
 matrix
function with nonnegative Hermitian values. Then
\begin{equation} \ln{\Bigg[ \det{\Bigg( \int\limits_{\Omega} M \ dP \Bigg)} \Bigg]} \geq \int\limits_{\Omega} \ln{[\det{M}]} \ dP . \label{3.18}
\end{equation}}

PROOF OF PART (ii) OF THEOREM \ref{tu32}.  For $k \in {\Bbb N}$ we define
$v_k : {\Bbb T} \to [1,\infty)$ via the rule
\begin{equation} v_k(t) := \left\{ \begin{array}{ll}
\| \underline{E}(t) \cdot \underline{F_k}(t) \| \ &, \mbox{ if} \quad \| \underline{E}(t) \cdot \underline{F_k}(t) \| \geq 1 \\
\hspace{0.2cm}&\hspace{0.2cm}\\
 \ &, \mbox{ if} \quad \| \underline{E}(t) \cdot \underline{F_k}(t) \| < 1 .
\end{array} \right. \label{3.19} \end{equation}
For $k \in {\Bbb N}$ and $t \in {\Bbb T}$ we then  have
\begin{equation} \ln{[v_k(t)]} \in [0,\infty) . \label{3.20}
\end{equation}
Combining $(\alpha)$ and (\ref{3.19}) we infer that for almost all
 $t \in {\Bbb T}$ with respect to $m$,
\begin{equation} \lim\limits_{k \to \infty} \ln{[v_k(t)]} = 0 . \label{3.21}
\end{equation}
For $k \in {\Bbb N}$ and $t \in {\Bbb T}$ we get the inequality
\[ \ln{[v_k(t)]} \leq \ln^+{\| \underline{E}(t) \|} +
\ \ln^+{\| \underline{F_k}(t) \|} \]
 from (\ref{3.19}),
which together with $(\beta)$ implies that the family $(\ln{v_k})_{k \in {\Bbb N}}$
is uniformly $m$-integrable. Combining this fact with (\ref{3.20}) and
(\ref{3.21}), an application of Vitali's Theorem provides
\begin{equation} \lim\limits_{k \to \infty} \int\limits_{{\Bbb T}} \ln{[v_k(t)]} \ m(dt) = 0 . \label{3.22}
\end{equation}
For $k \in {\Bbb N}$ we define $\Psi_k : {\Bbb D} \to {\Bbb C}$ via the
formula
\begin{equation} \Psi_k(z) := \exp \left\{ - \int\limits_{{\Bbb T}} \ln{[v_k(t)]} \ \frac{t+z}{t-z} \ m(dt) \right\} . \label{3.23}
\end{equation}
Therefore, in view of  (\ref{3.20}),  we obtain the inequality
\begin{eqnarray} |\Psi_k(z)| & = & \exp \left\{ \Re \left[ - \int\limits_{{\Bbb T}} \ln{[v_k(t)]} \ \frac{t+z}{t-z} \ m(dt) \right] \right\} \nonumber\\
& = & \exp \left\{ - \int\limits_{{\Bbb T}} \ln{[v_k(t)]} \ \frac{1-|z|^2}{|t-z|^2} \ m(dt) \right\} \leq \exp{\{0\}} = 1   \label{3.24}
\end{eqnarray}
for $z \in {\Bbb D}$.
In view of (\ref{3.21}), an application of Lebesgue's dominated convergence
 theorem yields
\begin{equation} \lim\limits_{k \to \infty} \Psi_k(z) = 1 \label{3.25}
\end{equation}
for all $z \in {\Bbb D}$. For almost all $t \in {\Bbb T}$ with respect to $m$
we get from (\ref{3.24})
\begin{equation} |\underline{\Psi_k}(t)| \leq 1 \label{3.26}
\end{equation}
and hence, upon taking into account that formula (\ref{3.23}) implies that
\begin{equation} |\underline{\Psi_k}(t)| = [v_k(t)]^{-1}, \label{3.27}
\end{equation}
we see from (\ref{3.19}) and $(\alpha)$ that
\begin{equation} \lim\limits_{k \to \infty} |\underline{\Psi_k}(t)| = 1 . \label{3.28}
\end{equation}
For $k \in {\Bbb N}$,
\begin{equation} \int\limits_{{\Bbb T}} |\underline{\Psi_k}(t)-1|^2 \ m(dt) = \int\limits_{{\Bbb T}} |\underline{\Psi_k}(t)|^2 \ m(dt) - 2 \Re{[\Psi_k(0)]} + 1 . \label{3.29}
\end{equation}
In view of (\ref{3.26}) and (\ref{3.28}), Lebesgue's dominated convergence
theorem yields
\begin{equation} \lim\limits_{k \to \infty} \int\limits_{{\Bbb T}} |\underline{\Psi_k}(t)|^2 \ m(dt) = m({\Bbb T}) = 1 . \label{3.30}
\end{equation}
Combining (\ref{3.25}), (\ref{3.29}) and (\ref{3.30}) we obtain
\begin{equation} \lim\limits_{k \to \infty} \int\limits_{{\Bbb T}} |\underline{\Psi_k}(t)-1|^2 \ m(dt) = 0 . \label{3.31}
\end{equation}
In view of (\ref{3.31}), the F. Riesz - Fischer theorem provides a subsequence
$(\underline{\Psi_{l_k}})_{k \in {\Bbb N}}$ of $(\underline{\Psi_k})_{k
 \in {\Bbb N}}$ such that
\begin{equation}
 \lim\limits_{k \to \infty} \underline{\Psi_{l_k}}(t) = 1 \label{3.32}
\end{equation}
for almost all $t \in {\Bbb T}$ with respect to $m$. Suppose that
 $k \in {\Bbb N}$ and  define
\begin{equation} \Phi_k := E \cdot F_k \cdot \Psi_k . \label{3.33}
\end{equation}
Then, since $E \in {\goth N}_n^+({\Bbb D})$, $F_k \in H_n^{\infty}({\Bbb D})$
 and (\ref{3.24}) holds,
we get
\begin{equation} \Phi_k \in {\goth N}_n^+({\Bbb D}) . \label{3.34}
\end{equation}
For almost all $t \in {\Bbb T}$ with respect to $m$ it follows from
 (\ref{3.33}), (\ref{3.19}) and (\ref{3.27}) that
\begin{equation} \| \underline{\Phi_k}(t) \| = |\underline{\Psi_k}(t)| \cdot \| \underline{E}(t) \cdot \underline{F_k}(t) \| \leq |\underline{\Psi_k}(t)|
 \cdot v_k(t) = 1 . \label{3.35}
\end{equation}
Therefore the maximum principle of V.I. Smirnov
implies that
\begin{equation} \| \underline{\Phi_k}(z) \| \leq 1 \label{3.36}
\end{equation}
for all $z \in {\Bbb D}$. In particular,
\begin{equation} \Phi_k \in H_n^{\infty}({\Bbb D}) . \label{3.37}
\end{equation}
From (\ref{3.34}) and (\ref{3.35}) it follows that
\begin{equation} \underline{\Phi_k}^{\star}(t) \cdot \underline{\Phi_k}(t) \leq I_n \label{3.38}
\end{equation}
for almost all $t \in {\Bbb T}$ with respect to $m$. Combining (\ref{3.33}),
 $(\alpha)$ and (\ref{3.28}) we get
\begin{equation} \lim\limits_{k \to \infty} \underline{\Phi_k}^{\star}(t) \cdot \underline{\Phi_k}(t)
= \lim\limits_{k \to \infty} |\underline{\Psi_k}(t)|^2 [ \underline{E}(t) \underline{F_k}(t)]^{\star} [\underline{E}(t) \underline{F_k}(t)] = I_n .
 \label{3.39}
\end{equation}
From (\ref{3.32}), (\ref{3.33}) and $(\alpha)$  we now obtain
\begin{equation} \lim\limits_{k \to \infty} \underline{\Phi_{l_k}}(t)
= \lim\limits_{k \to \infty} \underline{E}(t) \cdot \underline{F_{l_k}}(t) \cdot \underline{\Psi_{l_k}}(t) = I_n . \label{3.40}
\end{equation}
Using (\ref{3.37}), (\ref{3.38}), (\ref{3.40}) and Lebesgue's dominated convergence
theorem we get
\begin{equation} \lim\limits_{k \to \infty} \Phi_{l_k}(0) = \lim\limits_{k \to \infty} \int\limits_{{\Bbb T}} \underline{\Phi_{l_k}}(t) \ m(dt) = I_n . \label{3.41}
\end{equation}
Suppose that $k \in {\Bbb N}$. We define $M_k : {\Bbb T} \to {\goth M}_n$ via
the rule
\begin{equation} M_k(t) := \underline{\Phi_k}^{\star}(t) \cdot \underline{\Phi_k}(t) . \label{3.42}
\end{equation}
Then (\ref{3.42}) and (\ref{3.38}) imply that the inequality $0 \leq M_k(t) \leq I_n$
holds true for almost all $t \in {\Bbb T}$ with respect to $m$. Hence,
\begin{equation} 0 \leq \int\limits_{{\Bbb T}} M_k(t) \ m(dt) \leq I_n . \label{3.43}
\end{equation}
Now we apply the Generalized Minkowski inequality to the $M_k$. (Note that
Lebesgue measure $m$ is a probability measure.) From (\ref{3.43}) we infer
 first that
\begin{equation} \ln{\Bigg[ \det{\Bigg( \int\limits_{{\Bbb T}} M_k(t) \ m(dt) \Bigg)} \Bigg]} \leq \ln{[\det{I_n}]} = 0 . \label{3.44}
\end{equation}
Hence, (\ref{3.44}) and the Generalized Minkowski inequality guarantee that
\begin{equation} \int\limits_{{\Bbb T}} \ln{(\det{[M_k(t)]})} \ m(dt) \leq 0 . \label{3.45}
\end{equation}
Using (\ref{3.42}) and (\ref{3.33}) it follows that
\begin{eqnarray} \frac{1}{2} \ln{(\det{[M_k(t)]})} & =
 & \ln{|\det{[\underline{\Phi_k}(t)]}|} \nonumber\\
& = & \ln{|\det{[\underline{E}(t)]}|} + \ln{|\det\{ [\underline{F_k}(t)] \cdot [\underline{\Psi_k}(t)] \} |} \label{3.46}
\end{eqnarray}
for almost all $t \in {\Bbb T}$ with respect to $m$. Thus, from (\ref{3.45})
 and (\ref{3.46})
 we see that
\begin{equation} \int\limits_{{\Bbb T}} \ln{|\det{[\underline{E}(t)]}|} \ m(dt) \leq - \int\limits_{{\Bbb T}} \ln{|\det\{ [\underline{F_k}(t)] \cdot [\underline{\Psi_k}(t)] \} |} \ m(dt) . \label{3.47}
\end{equation}
By assumption, $F_k \in H_n^{\infty}({\Bbb D})$. Using (\ref{3.23}) and
 (\ref{3.24})
we see that $\Psi_k \in H_n^{\infty}({\Bbb D})$. Thus,\\ $F_k \cdot \Psi_k \in H_n^{\infty}({\Bbb D})$
and, consequently, $\det{[F_k \cdot \Psi_k]} \in H^{\infty}({\Bbb D})$.
Now Jensen's inequality gives
\begin{equation} - \int\limits_{{\Bbb T}} \ln{|\det\{ [\underline{F_k}(t)] \cdot [\underline{\Psi_k}(t)] \} |} \ m(dt) \leq - \ln{|\det\{ [\underline{F_k}(0)] \cdot [\underline{\Psi_k}(0)] \} |} . \label{3.48}
\end{equation}
From (\ref{3.47}) and (\ref{3.48})  it now  follows that
\begin{equation} \int\limits_{{\Bbb T}} \ln{|\det{[\underline{E}(t)]}|} \ m(dt) \leq - \ln{|\det\{ [\underline{F_k}(0)] \cdot [\underline{\Psi_k}(0)] \} |} . \label{3.49}
\end{equation}
From (\ref{3.33}) and (\ref{3.41}) we obtain
\begin{equation} \lim\limits_{k \to \infty} \ln{|\det\{ [F_{l_k}(0)] \cdot [\Psi_{l_k}(0)] \} |} = - \ln{| \det{[E(0)]} |} . \label{3.50}
\end{equation}
Combining (\ref{3.49}) and (\ref{3.50}) we obtain
\[ \int\limits_{{\Bbb T}} \ln{|\det{[\underline{E}(t)]}|} \ m(dt) \leq \ln{|\det{[E(0)]}|} . \]
By assumption, $E \in {\goth N}_n^+({\Bbb D})$. Thus, $\det{E} \in {\goth N}^+({\Bbb D})$
and Jensen's inequality yields
\[ \ln{|\det{[E(0)]}|} \leq \int\limits_{{\Bbb T}} \ln{| \det{[\underline{E}(t)]}|} \ m(dt) . \]
Hence,
\begin{equation} \int\limits_{{\Bbb T}} \ln{|\det{[\underline{E}(t)]}|} \ m(dt) = \ln{|\det{[E(0)]}|} . \label{3.51}
\end{equation}
From (\ref{3.51}) and Lemma \ref{lu21} we see that $\det{E} \in {\goth E}({\Bbb D})$.
Therefore, by definition \ref{du31}, $E \in {\goth E}_n({\Bbb D})$.
Part (ii) of Theorem \ref{tu32} is now  proved. \hfill $\boxx$\\[-0.3cm]

\begin{theo}
\label{tu32a}
{\rm (i) } {\sl Let $E \in {\goth E}_n({\Bbb D})$. Then there exists a sequence
$(F_k)_{k \in {\Bbb N}}$ from $H_n^{\infty}({\Bbb D})$ with the following properties:\\[0.3cm]
\hspace*{1.0cm}$(\alpha)$ For almost all $t \in {\Bbb T}$ with respect to $m$,
$
 \lim\limits_{k \to \infty} \underline{F_k}(t) \cdot \underline{E}(t) = I_n .
 $\\[0.2cm]
\hspace*{1.0cm}$(\beta)$ The family $\left(\ln^+{\|\underline{F_k}\|}\right)_{k \in {\Bbb N}}$ is
uniformly integrable with respect to {\rm m}.\\[0.2cm]
\hspace*{1.0cm}$(\gamma)$\parbox[t]{13.0cm}{ There exists a Borel subset $B_0$ of ${\Bbb T}$
 with $m(B_0)=0$
such that for all $k \in {\Bbb N}$ and all $t \in {\Bbb T} \setminus B_0$ the
 inequality
$
 \| \underline{F_k}(t) \cdot \underline{E}(t) \| \leq 1
$
holds.}\\[0.3cm]
\hspace*{0.6cm}{\rm (ii) } Let $E \in {\goth N}_n^+({\Bbb D})$ be such that
 there exists a sequence $(F_k)_{k \in {\Bbb N}}$
which belongs  to\\
\hspace*{1.2cm} $H_n^{\infty}({\Bbb D})$ and  satisfies
the above conditions $(\alpha)$ and $(\beta)$. Then $E \in {\goth E}_n({\Bbb D})$.
}
\end{theo}\\[-0.3cm]

 PROOF.  Combine Theorem \ref{tu32} and Remark \ref{ru31a}. \hfill $\boxx$

It should be mentioned that Ginzburg [Gi1] obtained a multiplicative
integral representation for outer functions which belong to
 ${\goth E}_n({\Bbb D})$.\\[0.2cm]

%%%%%%%%%%%%%%%%%%%%%%%%%%%%%%%%%%%%%%%%%%%%%%%%%%%%%%%%%%%%%%%%%%%%%%%%%%%%%%%
%%%%%%%%%%%%%%%%%%%%%%%%%%%%% SECTION 4 %%%%%%%%%%%%%%%%%%%%%%%%%%%%%%%%%%%%%%
%%%%%%%%%%%%%%%%%%%%%%%%%%%%%%%%%%%%%%%%%%%%%%%%%%%%%%%%%%%%%%%%%%%%%%%%%%%%%%%

\begin{minipage}{15.0cm}
\section{\hspace{-0.4cm}. MATRIX-VALUED INNER FUNCTIONS}
\end{minipage}\\[-0.1cm]
\setcounter{equation}{0}

In this section, we draw our attention to a distinguished subclass of the
Schur class ${\goth S}_{n \times n}$ (compare Definition \ref{du21}).\\[-0.2cm]

\begin{defi}
\label{du41}
Let $\Theta \in {\goth S}_{n \times n}({\Bbb D})$. Then $\Theta$ is called
{\sl inner} if
\begin{equation} I_n - \underline{\Theta}^{\star}(t) \cdot
 \underline{\Theta}(t)
 = {\Bbb O}_{n \times n} \label{4.1}
\end{equation}
for almost all $t \in {\Bbb T}$ with respect to $m$. The class of all
 $n \times n$ matrix-valued
inner functions will be denoted by ${\goth I}_n({\Bbb D})$.
\end{defi}\\[-0.2cm]

\begin{rem}
\label{ru41}
Let $\Theta \in {\goth I}_n({\Bbb D})$. Then obviously $\det{\Theta}
 \not\equiv 0$.
\end{rem}\\[-0.2cm]

\begin{rem}
\label{ru41a}
Let $\Theta \in {\goth I}_n({\Bbb D})$. Then $\Theta^{\top} \in {\goth I}_n({\Bbb D})$.
\end{rem}\\[-0.2cm]

The class ${\goth I}_n({\Bbb D})$ contains two important subclasses, namely
the so - called singular inner functions and the Blaschke-Potapov products.
Now we will formulate the corresponding definitions.

\begin{defi}
\label{du42}
Let $S \in {\goth I}_n({\Bbb D})$. Then $S$ is called  {\sl singular}, if
$\det{[S(z)]} \not= 0$ for all $z \in {\Bbb D}$ (or in other words if $S^{-1}$
is holomorphic  in ${\Bbb D}$). The class of all $n \times n$ matrix-valued
singular inner functions will be denoted by ${\goth I}_{n,s}({\Bbb D})$.
\end{defi}\\[-0.2cm]

\begin{rem}
\label{ru42}
If $S \in {\goth I}_{n,s}({\Bbb D})$, then $S^{-1} \in {\goth N}_n({\Bbb D})$,
be cause $S^{-1}$ admits the representation $S^{-1} = L \cdot (\det{S})^{-1}$
 with bounded holomorphic functions $L$ and $\det{S}$.
\end{rem}\\[-0.2cm]

\begin{lemma}
\label{lu41}
{\sl
Let $S \in {\goth I}_{n,s}({\Bbb D})$ be such that $S^{-1} \in
 {\goth N}_n^+({\Bbb D})$.
Then $S$ is constant.
}
\end{lemma}\\[-0.2cm]

 PROOF.  Since
 $\underline{S}(t)$ is unitary for a.e. $t \in {\Bbb T}$
it follows that
\[ \| \underline{S^{-1}}(t) \| = 1 . \]
Therefore, by the maximum principle of V.I. Smirnov,
 $\|S^{-1}(z)\| \leq 1$ for all
$z \in {\Bbb D}$. Since $\|S(z)\| \leq 1$ then it follows that $S(z)$ is a
unitary matrix for all $z \in {\Bbb D}$. However a holomorphic matrix function
with unitary values is necessarily constant (see e.g. Corollary 2.3.2
in [DFK]). \mbox{\ }  \hfill $\boxx$

Now we are going to define Blaschke-Potapov products. For this reason, we
 recall first
the notion of a scalar elementary Blaschke factor. Let $a \in {\Bbb D}$. Then
we define $b_a : {\Bbb D} \to {\Bbb C}$ via the rule
\begin{equation} b_a(z) := \left\{ \begin{array}{ll}
\frac{\displaystyle \raisebox{0.05cm}{$|a|$}}
{\displaystyle\raisebox{-0.05cm}{$a$}} \cdot \frac{
\displaystyle \raisebox{0.05cm}{$ a-z$}}
{\displaystyle\raisebox{-0.08cm} {$1-\overline{a}z$}}
 & , \mbox{ if } \quad a \in {\Bbb D} \setminus \{0\} \\
\hspace{0.2cm}z & , \mbox{ if } \quad a=0 \end{array} \right. . \label{4.2}
\end{equation}
Assume that $P \in {\goth M}_n$ is a non-zero orthoprojection matrix, i.e.,
that the conditions
\begin{equation} P^2=P \qquad \qquad P=P^{\star} \label{4.3}
\end{equation}
are satisfied. Then the matrix-valued function $B_{a,P} : {\Bbb D} \to
 {\goth M}_n$
which is defined by
\begin{equation} B_{a,P}(z) := I_n + [b_a(z)-1] \cdot P \label{4.4}
\end{equation}
is called {\em the Blaschke-Potapov elementary factor associated with $a$
 and $P$.}\\
From (\ref{4.3}) and (\ref{4.4}) it is clear that
\begin{equation} \det{[B_{a,P}]} = (b_a)^{\mbox{{\scriptsize rank }} P} .
\label{4.5}
\end{equation}
Suppose that $(z_k)_{k \in I}$ is a sequence from ${\Bbb D}$ and that $(P_k)_{k \in I}$
is a sequence of orthoprojection matrices for which the condition
\begin{equation}  \sum\limits_{k \in I} (1-|z_k|) \cdot \mbox{tr } P_k <
 +\infty \label{4.6}
\end{equation}
is fulfilled. (The index set $I$ can be finite or infinite.)
Then, according to a result due to V.P. Potapov [Pot],   the product
\begin{equation} \prod\limits_{k \in I}^{\displaystyle\curvearrowright} B_{z_k,P_k}(z)
 \qquad \qquad
\left( \mbox{resp. } \ \ \prod\limits_{k \in I}^{\displaystyle\curvearrowleft}
 B_{z_k,P_k}(z) \right) \label{4.7}
\end{equation}
converges for all $z \in {\Bbb D}$. (V.P. Potapov has also shown that
 condition (\ref{4.6}) is necessary for the convergence of the product in
 (\ref{4.7})).

\begin{defi}
\label{du43}
{\sl
Let $B: {\Bbb D} \to {\goth M}_n$. Then $B$ is called a {\rm left} (resp. {\rm right})
{\rm Blaschke-Potapov product} if $B$ is a constant function with unitary
value or if there exist a unitary matrix $V$, a set of orthoprojection
 matrices $(P_k)_{k \in I}$ and sequences $(z_k)_{k \in I}$
which belong to ${\Bbb D}$  such that
(\ref{4.6}) is satisfied and moreover the representation
\[ B=\left( \prod\limits_{k \in I}^{\displaystyle\curvearrowright} B_{z_k,P_k}(z) \right) \cdot V \qquad \left( \mbox{resp. } \ \ B=V \cdot \left( \prod\limits_{k \in I}^{\displaystyle\curvearrowleft} B_{z_k,P_k}(z) \right) \right) \]
is valid. The set of left (resp. right) Blaschke-Potapov products will be denoted
by ${\goth I}_{n,B,l}({\Bbb D})$ (resp. ${\goth I}_{n,B,r}({\Bbb D})$).
}
\end{defi}\\[-0.2cm]

We will see below that each left Blaschke-Potapov product is also a right
Blaschke-Potapov product and vice versa. Moreover, it will turn out that
${\goth I}_{n,B,l}({\Bbb D}) \subseteq {\goth I}_n({\Bbb D})$.\\[-0.2cm]

\begin{lemma}
\label{lu42}
{\sl
Let $A,B \in {\goth C}_n$ be such that $A \cdot B$ is unitary.
Then $A$ and $B$ are unitary too.}
\end{lemma}

 PROOF.  Since  $A,B \in {\goth C}_n$ we have $I_n-AA^{\star}
 \geq {\Bbb O}_{n \times n}$
and $I_n-BB^{\star} \geq {\Bbb O}_{n \times n}$. Hence, $A(I_n-BB^{\star})A^{\star} \geq {\Bbb O}_{n \times n}$.
Therefore, the identity
\[ {\Bbb O}_{n \times n} = I_n - (AB)(AB)^{\star} = (I_n -AA^{\star})+A(I_n-BB^{\star})A^{\star} \]
implies that $I_n-AA^{\star} = {\Bbb O}_{n \times n}$ and $A(I_n-BB^{\star})A^{\star} = {\Bbb O}_{n \times n}$.
Thus, $A$ is unitary. In particular, we have $\det{A} \not= 0$.
This implies that $I_n-BB^{\star} = {\Bbb O}_{n \times n}$ and hence that
 $B$ is unitary too. \qquad \hfill $\boxx$\\[0.4cm]

\begin{theo}
\label{tu41}
{\sl
Suppose that $\Theta \in {\goth I}_n({\Bbb D})$.\\[-0.5cm]

{\rm (a) }\parbox[t]{14.6cm}{There exist functions $B \in {\goth I}_{n,B,l}({\Bbb D})$ (resp. $C \in {\goth I}_{n,B,r}({\Bbb D})$)
and $S \in {\goth I}_{n,s}({\Bbb D})$\\ (resp. $T \in {\goth I}_{n,s}({\Bbb D})$) such that
the multiplicative representation
\begin{equation} \Theta = B \cdot S \qquad \qquad (\mbox{resp. } \ \Theta =T \cdot C) \label{4.8}
\end{equation}
holds true.}\\[-0.2cm]

{\rm (b) }\parbox[t]{14.6cm}{ Suppose that the functions $B_1,B_2 \in {\goth I}_{n,B,l}({\Bbb D})$
(resp. $C_1,C_2 \in {\goth I}_{n,B,r}({\Bbb D})$) and $S_1,S_2 \in {\goth I}_{n,s}({\Bbb D})$
(resp. $T_1,T_2 \in {\goth I}_{n,s}({\Bbb D})$) satisfy $B_1S_1=B_2S_2=\Theta$
(resp. $T_1C_1=T_2C_2=\Theta$). Then there exist a unitary matrix $U \in {\goth M}_n$
(resp. $V \in {\goth M}_n$) such that $B_2=B_1U$ and $S_2=U^{\star}S_1$
(resp. $C_2=VC_1$ and $T_2=T_1V^{\star}$) are fulfilled.}
}
\end{theo}\\[-0.1cm]

 PROOF.  Theorem \ref{tu41} is a special case of a much more general result
due to V.P. Potapov [Pot]. The Potapov theory handles the case of meromorphic
matrix-valued functions in ${\Bbb D}$ which have a nonidentically vanishing
 determinant
and which are $J$-contractive where $J$ is a signature matrix (i.e. $J=J^{\star}$
and $J^2=I_n$). In the special case that $J=I$, V.P. Potapov's result
 (see [Pot] and also
a series of papers by Ginzburg [Gi1] - [Gi5], [GiSh]) provides the existence of functions $B \in {\goth I}_{n,B,l}({\Bbb D})$
and $S \in {\goth S}_{n \times n}({\Bbb D})$ such that
\begin{equation} \Theta = B \cdot S \label{4.9}
\end{equation}
and for all $z \in {\Bbb D}$,
\begin{equation} \det{[S(z)]} \not= 0 . \label{4.10}
\end{equation}
Since the boundary function $\underline{\Theta}$ has unitary values almost
everywhere with respect to $m$ we infer from Lemma \ref{lu42}  that the
 boundary functions $\underline{B}$
and $\underline{S}$ also  have unitary values almost everywhere with respect
to $m$. Taking
 into
account (\ref{4.10}) we obtain $S \in {\goth I}_{n,s}({\Bbb D})$. The uniqueness part
goes back to V.P. Potapov [Pot] too. \qquad \hfill $\boxx$\\[-0.2cm]

\begin{lemma}
\label{lu43}
{\sl
Let $M \in {\goth C}_n$. Then\\[0.1cm]
\hspace*{0.6cm}{\rm (a) } $|\det{M}| \leq 1$\\[0.1cm]
\hspace*{0.6cm}{\rm (b) } $|\det{M}| = 1$ if and only if $M$ is unitary.
}
\end{lemma}\\[-0.2cm]

 PROOF. Let  $(l_k(M^{\star}M))_{k=1}^n$ denote  the system of eigenvalues
of $M^{\star}M$. Then, since  $M \in {\goth C}_n$,
 $0 \leq l_k(M^{\star}M) \leq 1$
for all $k \in \{ 1, \ldots , n \}$. Thus, as
\begin{equation} |\det{M}|^2 = \det{(M^{\star}M)} =
 \prod\limits_{k=1}^n l_k(M^{\star}M), \label{4.11}
\end{equation}
we see that $|\det{M}| \leq 1$ with equality  if and only if
\ $l_k(M^{\star}M) = 1$ for all $k \in \{ 1, \ldots , n \}$.
But $l_k(M^{\star}M)=1$ for all $k \in \{ 1, \ldots , n \}$ if and only
if $M^{\star}M=I_n$. \qquad \hfill $\boxx$

Now we recall a well-known characterization of Blaschke products (see e.g.,
Privalov [Pri, Ch.I, Sec.7.1]).

\begin{lemma}
\label{lu44}
{\sl Let $\Theta \in {\goth S}_{1 \times 1}({\Bbb D})$. Then $\Theta$ is a
 Blaschke product
if and only if}
\[ \lim\limits_{r \to 1-0} \int\limits_{{\Bbb T}} \ln{|\det{[\Theta(rt)]}|} \ m(dt) = 0 . \]
\end{lemma}

\begin{theo}
\label{tu42}
{\sl
Let $f \in {\goth S}_{n \times n}({\Bbb D})$. Then:\\[0.1cm]
\hspace*{0.8cm}{\rm (a) } The function $\det{f}$ belongs to ${\goth S}_{1 \times 1}({\Bbb D})$.\\
\hspace*{0.8cm}{\rm (b) } $f \in {\goth I}_n({\Bbb D})$ if and only if $\det{f} \in {\goth I}_1({\Bbb D})$.
If $f \in {\goth I}_n({\Bbb D})$ then $\det{f} \not\equiv 0$.\\
\hspace*{0.8cm}{\rm (c) } $f \in {\goth I}_{n,s}({\Bbb D})$ if and only if $\det{f} \in {\goth I}_{1,s}({\Bbb D})$.\\
\hspace*{0.8cm}{\rm (d) } The following statements are equivalent:\\
\hspace*{1,4cm}
\parbox[t]{13.0cm}
{
 {\rm (i) } $f \in {\goth I}_{n,B,l}({\Bbb D})$,\\
                {\rm (ii) } $f \in {\goth I}_{n,B,r}({\Bbb D})$,\\
                {\rm (iii) } $\det{f}$ is a Blaschke product.\\
                {\rm (iv) } The limit relation
$\displaystyle
 \lim\limits_{r \to 1-0} \int\limits_{{\Bbb T}} \ln{|\det{[f(rt)]}|} \ m(dt) = 0
$
holds true.}
}
\end{theo}

 PROOF.  The assertions stated in part (a) and (b) are immediate consequences
 of Lemma \ref{lu43}.
Part (c) follows from part (a) and the definition of a singular inner function.\\
It remains to prove part (d). From (a) and Lemma \ref{lu44} we can immediately
conclude the equivalence of statements (iii) and (iv). In view of (\ref{4.5}),
 it is readily
checked that each of the conditions (i) and (ii) implies (iii). Now suppose
that (iii) holds. By virtue of part (b) we see that $f$ is an inner function.
>From Theorem \ref{tu41} we infer that there exist functions $B \in {\goth I}_{n,B,l}({\Bbb D})$
and $S \in {\goth I}_{n,s}({\Bbb D})$ satisfying the multiplicative decomposition
$f=B \cdot S$. Hence, $\det{f}=\det{B} \cdot \det{S}$. The implication
``(i) $\Rightarrow$ (iii)'' which is already verified shows that $\det{B}$ is a
Blaschke product. Part (c) yields that $\det{S}$ is a singular inner function.
Therefore, the uniqueness part of Theorem \ref{tu41} yields that $\det{S}$ is a
constant inner function with unimodular value. Hence, we obtain from part (b)
of Lemma \ref{lu43} that  the matrix $S(z)$ is unitary for each
 $z \in {\Bbb D}$. Since $S$ belongs to ${\goth S}_{n \times n}({\Bbb D})$,
 the maximum
modulus principle for matrix-valued Schur functions (see e.g. [DFK,Corollary 2.3.2])
implies that $S$ is a constant function. From $f=B \cdot S$ we infer that
(i) holds. The implication ``(iii) $\Rightarrow$ (ii)'' can be shown analogously.
The theorem is proved. \qquad \hfill $\boxx$

For further results on matrix-valued and operator-valued inner functions
we refer the reader to the monographs Helson [Hel1], Sz.-Nagy and Foias
[SZNF] and Nikolskii [Nik2].\\[0.2cm]

%%%%%%%%%%%%%%%%%%%%%%%%%%%%%%%%%%%%%%%%%%%%%%%%%%%%%%%%%%%%%%%%%%%%%%%%%%%%%%%
%%%%%%%%%%%%%%%%%%%%%%%%%% SECTION 5 %%%%%%%%%%%%%%%%%%%%%%%%%%%%%%%%%%%%%%%%%%
%%%%%%%%%%%%%%%%%%%%%%%%%%%%%%%%%%%%%%%%%%%%%%%%%%%%%%%%%%%%%%%%%%%%%%%%%%%%%%%
\begin{minipage}{15.0cm}
\section{\hspace{-0.4cm}. INNER - OUTER FACTORIZATION}
\end{minipage}\\[-0.1cm]
\setcounter{equation}{0}

This section is aimed at a Smirnov class generalization of the inner-outer
factorization of matrix-valued functions belonging to the Hardy class
$H_n^2({\Bbb D})$.

Let us recall the following notions:

\begin{defi}
\label{du51}
{\sl
{\rm The Hardy class $H_n^2({\Bbb D})$}  is the set of all matrix-valued functions
$F: {\Bbb D} \to {\goth M}_n$ which are holomorphic in ${\Bbb D}$ and satisfy
\[ \sup\limits_{r \in [0,1)} \int\limits_{{\Bbb T}} \| F(rt) \|^2 \ m(dt) < \infty . \]
}
\end{defi}

\begin{rem}
\label{ru51}
Obviously,\ \ $H_n^{\infty}({\Bbb D}) \subseteq H_n^2({\Bbb D}) \subseteq
 {\goth N}_n^+({\Bbb D})$.
\end{rem}

\begin{rem}
\label{ru52}
Define $\| \bullet \|_{H^2} : H_n^2({\Bbb D}) \to [0,\infty)$ via
\[ F \to \sqrt{\sup\limits_{r \in [0,1)} \int\limits_{{\Bbb T}} \| F(rt)
 \|^2 \ m(dt) } . \]
Then $(H_n^2({\Bbb D}), \| \bullet \|_{H^2})$ is a complex Hilbert space.
\end{rem}

\begin{rem}
\label{ru52aa}
Let $S \in {\goth S}_{n \times n}({\Bbb D})$ be such that $\det{(I_n+S)}$
does not identically vanish in ${\Bbb D}$. Then $(I_n+S)
 \in {\goth E}_n({\Bbb D}) \cap H_n^{\infty}({\Bbb D})$
(see Arov [Ar1], Lemma 3.1).
\end{rem}

The definition of a matrix-valued outer function given above
(see Definition \ref{du31}) is too rough for the purposes of prediction
theory of stationary sequences. For this reason, P.R. Masani [Ma1, Ma2]
introduced  the following notion for the space $H_n^2({\Bbb D})$  (compare
Lemma \ref{lu21a}).\\[-0.2cm]

\begin{defi}
\label{du52a}
{\sl
Let $E \in H_n^2({\Bbb D})$. Then $E$ is said to be {\rm left optimal}
 {\rm(} resp.
{{\rm right optimal}{\rm )}}
if $E$ has the following property:
If $F \in H_n^2({\Bbb D})$ satisfies $[\underline{F}(t)] \cdot
 [\underline{F}(t)]^* = [\underline{E}(t)] \cdot [\underline{E}(t)]^*$
(resp. $[\underline{F}(t)]^* \cdot [\underline{F}(t)] = [\underline{E}(t)]^*
 \cdot [\underline{E}(t)]$)
then $[F(0)] \cdot [F(0)]^* \leq [E(0)] \cdot [E(0)]^*$ (resp.
$[F(0)]^* \cdot [F(0)] \leq [E(0)]^* \cdot [E(0)]$).
}
\end{defi}\\[-0.2cm]

\begin{rem}
\label{ru52a}
Let $E \in H_n^2({\Bbb D})$. Then $E$ is left optimal if and only if
$E^{\top}$ is right optimal.
\end{rem}

This notion of optimality is closely related to the following definition
which in the scalar case goes back to Beurling [Be].\\[-0.2cm]

\begin{defi}
\label{du52}
{\sl
Let $E \in H_n^2({\Bbb D})$. Then $E$ is called {\rm left Beurling-outer}
{\rm(} resp. {\rm right Beurling outer} {\rm)} if there exists a sequence
 $(f_k)_{k \in {\Bbb N}}$ from $H_n^{\infty}({\Bbb D})$
which satisfies
\[ \lim\limits_{k \to \infty} \int\limits_{{\Bbb T}}
 \| \underline{F_k}(t) \cdot
 \underline{\raisebox{0.0cm}[0.0cm][0.06cm]{$E$}}(t) -
 I_n \|^2 \ m(dt)
 = 0 \quad
\mbox {\rm\Big(} \ resp. \
 \lim\limits_{k \to \infty} \int\limits_{{\Bbb T}} \|
 \underline{\raisebox{0.0cm}[0.0cm][0.06cm]{$E$}}(t)
 \cdot \underline{F_k}(t) - I_n \|^2 \ m(dt) = 0\Big ). \]
The class of all $n \times n$ matrix-valued left Beurling-outer (resp. right
Beurling outer) functions will be denoted by ${\goth E}_{n,B,l}({\Bbb D})$
(resp. ${\goth E}_{n,B,r}({\Bbb D})$).
}
\end{defi}\\[-0.2cm]

\begin{rem}
\label{ru52b}
Let $E \in H_n^2({\Bbb D})$. Then $E \in {\goth E}_{n,B,l}({\Bbb D})$
if and only if $E^{\top} \in {\goth E}_{n,B,r}({\Bbb D})$.
\end{rem}

\begin{rem}
\label{ru53}
Let $E \in H_n^2({\Bbb D})$. Then it is readily checked that $E$ is left
Beurling outer (resp. right Beurling outer) if and only if the subspace
$H_n^2({\Bbb D}) \cdot E$ (resp. $E \cdot H_n^2({\Bbb D})$) is dense in
$(H_n^2({\Bbb D}),\| \bullet \|_{H^2})$.
\end{rem}

\begin{rem}
\label{ru54}
Let $E$ be a function belonging to ${\goth E}_{n,B,l}({\Bbb D})$ or
${\goth E}_{n,B,r}({\Bbb D})$. Then for all $z \in {\Bbb D}$ the relation
$\det{[E(z)]} \not=0$ holds true.
\end{rem}

 PROOF :  Let us consider the case $E \in {\goth E}_{n,B,r}({\Bbb D})$.
Then there exists a sequence $(F_k)_{k \in {\Bbb N}}$ from
 $H^{\infty}({\Bbb D})$ such that
\[ \lim\limits_{k \to \infty} \int\limits_{{\Bbb T}} \|
 \underline{\raisebox{0.0cm}[0.0cm][0.06cm]{$E$}}(t) \cdot \underline{F_k}(t) - I_n \|^2 \ m(dt) = 0 . \]
From this it follows  by the Poisson integral
representation for $H_n^2({\Bbb D})$ functions that
\[ \lim\limits_{k \to \infty} E(z) \cdot F_k(z) = I_n \]
for $z \in {\Bbb D}$ and hence that
\[ \lim\limits_{k \to \infty} \det{[E(z)]} \cdot \det{[F_k(z)]} = 1. \]
Thus, $\det{[E(z)]} \not= 0$. If $E \in {\goth E}_{n,B,l}({\Bbb D})$,
then the assertion follows from Remark \ref{ru52b} and the  preceding
analysis. \hfill $\boxx$

The following result due to Masani [Ma2, Corollary 4.6] clarifies
the relation between optimality and Beurling-outerness.\\[-0.2cm]

\begin{theo}
\label{tu50}
{\sl
Let $E \in H_n^2({\Bbb D})$. Then:\\
\hspace*{0.8cm}{\rm (a) } \parbox[t]{14.0cm}{If $\det{E} \not\equiv 0$ and $E$
 is left optimal (resp.
right optimal), then $E \in {\goth E}_{n,B,l}({\Bbb D})$
(resp. $E \in {\goth E}_{n,B,r}({\Bbb D})$).}\\
\hspace*{0.8cm}{\rm (b) }\parbox[t]{14.0cm} {If $E
 \in {\goth E}_{n,B,l}({\Bbb D})$ (resp. $E \in {\goth E}_{n,B,r}({\Bbb D})$),
then $E$ is left optimal (resp. right optimal).}
}
\end{theo}

The notion of optimality is more general than the notion of Beurling - outer
because it allows  the functions in question to have identically vanishing
determinants.  In the theory of multivariate stationary
stochastic processes this corresponds to  the case of a singular prediction
 error matrix.

The following result plays a key role in the theory of holomorphic
 matrix-valued functions.

\begin{theo}
\label{tu51}
{\sl
Let $F \in H_n^2({\Bbb D})$ be such that $\det{F} \not\equiv 0$. Then:\\[0.2cm]
\hspace*{0.8cm}{\rm (i) }\  \parbox[t]{14.0cm}{ There exist functions
 $\Theta_r \in {\goth I}_n({\Bbb D})$ and
$E_r \in {\goth E}_{n,B,r}({\Bbb D})$ such that the multiplicative
 decomposition
\[ F = \Theta_r \cdot E_r \]
is satisfied.}\\[0.2cm]
\hspace*{0.8cm}{\rm (ii) }\  \parbox[t]{14.0cm}{ Suppose that the functions
 $\Theta_{r1},\Theta_{r2} \in {\goth I}_n({\Bbb D})$
and $E_{r1},E_{r2} \in {\goth E}_{n,B,r}({\Bbb D})$ satisfy
\[ \Theta_{r1} \cdot E_{r1} = \Theta_{r2} \cdot E_{r2} = F . \]
Then there exists a unitary matrix $V \in {\goth M}_n$ such that
 $\Theta_{r2}=\Theta_{r1} \cdot V$
and $E_{r2}=V^{\star} \cdot E_{r1}$ are fulfilled.}\\[0.2cm]
\hspace*{0.8cm}{\rm (iii) }\  \parbox[t]{14.0cm}{ There exist functions
 $\Theta_l \in {\goth I}_n({\Bbb D})$ and
$E_l \in {\goth E}_{n,B,l}({\Bbb D})$ such that the multiplicative
 decomposition
\[ F = E_l \cdot \Theta_l \]
is satisfied.}\\[0.2cm]
\hspace*{0.8cm}{\rm (iv) }\  \parbox[t]{14.0cm}{ Suppose that the functions
 $\Theta_{l1},\Theta_{l2} \in {\goth I}_n({\Bbb D})$
and $E_{l1},E_{l2} \in {\goth E}_{n,B,l}({\Bbb D})$ satisfy
\[ E_{l1} \cdot \Theta_{l1}=E_{l2} \cdot \Theta_{l2} = F . \]
Then there exists a unitary matrix $U \in {\goth M}_n$ such that
 $\Theta_{l2}=U \cdot \Theta_{l1}$
and $E_{l2}=E_{l1} \cdot U^{\star}$ are fulfilled.}\\[0.2cm]
}
\end{theo}

Theorem \ref{tu51} was proved independently by several authors (see Masani
 [Ma2, 4.3, 4.4], Helson and Lowdenslager [HL2, Theorem 15], Rozanov [Roz1,
 Theorem 5]).  The Beurling-Lax-Halmos Theorem  (see Beurling [Be], Lax [La],
 Halmos [Hal] and also Masani [Ma2, Theorem 3.8.])
which describes the structure of shift invariant left (resp. right)
submodules of $H_n^2({\Bbb D})$ lies at the heart of the proof.\\[-0.2cm]

\begin{theo}
\label{tu52}
{\sl
The identities
\[ {\goth E}_{n,B,l}({\Bbb D}) = {\goth E}_{n,B,r}({\Bbb D}) =
 {\goth E}_n({\Bbb D}) \cap H_n^2({\Bbb D}) \]
are valid.\\[0.1cm]
}
\end{theo}

 PROOF.  First we show that
\[ {\goth E}_{n,B,r}({\Bbb D}) = {\goth E}_n({\Bbb D}) \cap H_n^2({\Bbb D}) .
 \]
Our proof is based  mainly  on Theorem \ref{tu32}.\\
First assume that $E \in {\goth E}_n({\Bbb D}) \cap H_n^2({\Bbb D})$.
Then part (i) of Theorem \ref{tu32} guarantees the existence of a sequence
$(F_k)_{k \in {\Bbb N}}$ from $H_n^{\infty}({\Bbb D})$ with the properties
$(\alpha)$, $(\beta)$ and $(\gamma)$ formulated there. In view of property
$(\gamma)$, there exists a Borel subset $B_0$ of ${\Bbb T}$ with $m(B_0)=0$
 such that for all
$k \in {\Bbb N}$ and all $t \in {\Bbb T} \setminus B_0$ the inequality
\begin{equation} \| \underline{E}(t) \cdot \underline{F_k}(t) - I_n \| \leq \|
 \underline{E}(t) \cdot \underline{F_k}(t) \| + \| I_n \| \leq 2 \label{5.1}
\end{equation}
holds . In view of  $(\alpha)$ and (\ref{5.1}), an application of
Lebesgue's dominated convergence theorem yields
\[ \lim\limits_{k \to \infty} \int\limits_{{\Bbb T}} \| \underline{E}(t) \cdot
 \underline{F_k}(t) - I_n \|^2 \ m(dt) = 0 . \]
Thus, $E \in {\goth E}_{n,B}({\Bbb D})$. Hence, the inclusion
\begin{equation} {\goth E}_n({\Bbb D}) \cap H_n^2({\Bbb D}) \subseteq
 {\goth E}_{n,B}({\Bbb D}) \label{5.2}
\end{equation}
holds true.\\
Now assume that $E \in {\goth E}_{n,B}({\Bbb D})$. Then Definition \ref{du52}
implies that
\begin{equation} E \in H_n^2({\Bbb D}) . \label{5.3}
\end{equation}
We will show that $E$ satisfies the conditions $(\alpha)$ and $(\beta)$ in
Theorem \ref{tu32}. In view of Definition \ref{5.2} there exists a sequence
$(F_k)_{k \in {\Bbb N}}$ from $H_n^{\infty}({\Bbb D})$ for which
\begin{equation} \lim\limits_{k \to \infty} \int\limits_{{\Bbb T}} \| \underline{E}(t) \cdot \underline{F_k}(t) - I_n \|^2 \ m(dt) = 0 . \label{5.4}
\end{equation}
Obviously, for $k \in {\Bbb N}$ and $t \in {\Bbb T}$ the inequality
\begin{equation} 0 \leq \ln^+{\| \underline{E}(t) \cdot \underline{F_k}(t) \|} \leq \| \underline{E}(t) \cdot \underline{F_k}(t) - I_n \| \label{5.5}
\end{equation}
holds true. From (\ref{5.4}) and (\ref{5.5}) it then follows that
\[ \lim\limits_{k \to \infty} \int\limits_{{\Bbb T}} \ln^+{\| \underline{E}(t) \cdot \underline{F_k}(t)\|} \ m(dt) = 0 . \]
Hence, the family $(\ln^+{\| \underline{E} \cdot \underline{F_k} \|} )_{k \in {\Bbb N}}$
is uniformly $m$ - integrable.
In view of Remark \ref{ru54} we see that $\det{[E(z)]} \not= 0$
for all $z \in {\Bbb D}$.
Since $E \in H_n^2({\Bbb D}) \subseteq {\goth N}_n({\Bbb D})$
 we now obtain $E^{-1} \in {\goth N}_n({\Bbb D})$.
Hence, $\ln{\|\underline{E}^{-1}\|}=\ln{\|\underline{E^{-1}}\|}$
is $m$-integrable. Clearly, for $k \in {\Bbb N}$ and $t \in {\Bbb T}$
the inequality
\begin{equation} \ln^+{\| \underline{F_k}(t)\|} \leq \ln^+{\| \underline{E}(t) \cdot \underline{F_k}(t)\|} + \ln^+{\| [\underline{E}(t)]^{-1} \|} \label{5.6}
\end{equation}
holds true. Since the family $(\ln^+{\| \underline{E} \cdot \underline{F_k} \|})_{k \in {\Bbb N}}$
is uniformly $m$-integrable and since $\ln{\| \underline{E}^{-1} \|}$ is
 $m$-integrable it follows
from (\ref{5.6})  that the family
 $(\ln^+{\| \underline{F_k} \|})_{k \in {\Bbb N}}$
is uniformly $m$-integrable. Taking into account (\ref{5.4}), the Theorem of
F. Riesz - Fischer provides the existence of a subsequence $(F_{l_k})_{k \in {\Bbb N}}$
of $(F_k)_{k \in {\Bbb N}}$ such that
\[ \lim\limits_{k \to \infty} \underline{E}(t) \cdot \underline{F_{l_k}}(t) = I_n \]
for $m$-almost all $t \in {\Bbb T}$. Since the family $(\ln^+{\|F_{l_k}\|})_{k \in {\Bbb N}}$
is also uniformly $m$-integrable the conditions $(\alpha)$ and $(\beta)$
in Theorem \ref{tu32} are satisfied for the sequence $(F_{l_k})_{k \in {\Bbb N}}$.
Thus, part (ii) of Theorem \ref{tu32} implies that
\begin{equation} E \in {\goth E}_n({\Bbb D}) . \label{5.7}
\end{equation}
From (\ref{5.3}) and (\ref{5.7}) we obtain ${\goth E}_{n,B}({\Bbb D}) \subseteq {\goth E}_n({\Bbb D}) \cap H_n^2({\Bbb D})$.\\
An application of (\ref{5.2}) shows that
\begin{equation} {\goth E}_{n,B,r}({\Bbb D}) = {\goth E}_n({\Bbb D}) \cap H_n^2({\Bbb D}) . \label{5.7a}
\end{equation}
From (\ref{5.7a}) and Remarks \ref{ru31a} and \ref{ru52b} we then get
\[ {\goth E}_{n,B,l}({\Bbb D}) = {\goth E}_n({\Bbb D}) \cap H_n^2({\Bbb D}) . \]
Thus, the theorem is proved.  \hfill $\boxx$\\[0.1cm]

\begin{theo}
\label{tu53}
{\sl
{\rm (Inner - outer factorization in the Smirnov class ${\goth N}_n^+({\Bbb D})$ )}. \
Let $F \in {\goth N}_n^+({\Bbb D})$ be such that $\det{F} \not\equiv 0$. Then:\\
\hspace*{0.8cm}{\rm (i) }\  \parbox[t]{14.0cm}{ There exist functions $\Theta_r \in {\goth I}_n({\Bbb D})$ and
$E_r \in {\goth E}_n({\Bbb D})$ such that
\[ F = \Theta_r \cdot E_r . \]}\\[0.2cm]
\hspace*{0.8cm}{\rm (ii) }\  \parbox[t]{14.0cm}{ Suppose that the functions $\Theta_{r1}, \Theta_{r2} \in {\goth I}_n({\Bbb D})$ and
$E_{r1}, E_{r2} \in {\goth E}_n({\Bbb D})$ satisfy
\[ \Theta_{r1} \cdot E_{r1} = \Theta_{r2} \cdot E_{r2} = F . \]
Then there exists a unitary matrix $V \in {\goth M}_n$ such that $\Theta_{r2}=\Theta_{r1} \cdot V$
and $E_{r2}=V^{\star} \cdot E_{r1}$ .}\\[0.2cm]
\hspace*{0.8cm}{\rm (iii) }\  \parbox[t]{14.0cm}{ There exist functions $\Theta_l \in {\goth I}_n({\Bbb D})$ and
$E_l \in {\goth E}_n({\Bbb D})$ such that
\[ F = E_l \cdot \Theta_l . \]}\\[0.2cm]
\hspace*{0.8cm}{\rm (iv) }\  \parbox[t]{14.0cm}{ Suppose that the functions $\Theta_{l1}, \Theta_{l2} \in {\goth I}_n({\Bbb D})$ and
$E_{l1}, E_{l2} \in {\goth E}_n({\Bbb D})$ satisfy
\[ E_{l1} \cdot \Theta_{r1} = E_{l2} \cdot \Theta_{r2} = F . \]
Then there exists a unitary matrix $U \in {\goth M}_n$ such that $\Theta_{l2}=U \cdot \Theta_{l1}$
and $E_{l2}=E_{l1} \cdot U^{\star}$.}
}\\[0.1cm]
\end{theo}

 PROOF.  We derive these results from Theorem \ref{tu51}.\\
(i) In view of Lemma \ref{lu23} there exist functions $d \in {\goth E}({\Bbb D})$
and $\Phi \in {\goth S}_{n \times n}({\Bbb D})$ such that
\begin{equation} F = \frac{1}{d} \cdot \Phi . \label{5.8}
\end{equation}
Since $\det{F} \not\equiv 0$,  it follows from (\ref{5.8}) that
$\det{\Phi} \not\equiv 0$.
 Thus as ${\goth S}_{n \times n}({\Bbb D}) \subseteq H_n^2({\Bbb D})$,
Theorem \ref{tu51} ensures the existence of functions $\Theta_r \in {\goth I}_n({\Bbb D})$
and $E_{r,B} \in {\goth E}_{n,B}({\Bbb D})$ such that
\begin{equation} \Phi = \Theta_r \cdot E_{r,B} . \label{5.9}
\end{equation}
We set
\begin{equation} E:=d \cdot E_{r,B} . \label{5.10}
\end{equation}
According to Theorem \ref{tu52} it follows that $E_{r,B}
 \in {\goth E}_n({\Bbb D})$.
Since $d \in {\goth E}({\Bbb D})$ we get $E \in {\goth E}_n({\Bbb D})$ from
 (\ref{5.10}).
Thus (i) is proved.\\
(ii) The factorizations $F = \Theta_{r1} \cdot E_{r1} = \Theta_{r2} \cdot E_{r2}$
yield  the factorizations
\begin{equation} \Theta_{r1} \cdot E_{r1,B} = \Theta_{r2} \cdot E_{r2,B} = \Phi , \label{5.11}
\end{equation}
upon setting $E_{r1,B}:=d \cdot E_{r1}$ ,  $E_{r2,B}:=d \cdot E_{r2}$
and invoking (\ref{5.8}).\\
From $\Phi \in {\goth S}_{n \times n}({\Bbb D})$, (\ref{5.11}) and its
 definition it is clear that
\[ E_{r1,B} , E_{r2,B} \in {\goth E}_n({\Bbb D}) \cap {\goth S}_{n \times n}({\Bbb D}) . \]
Thus, from Theorem \ref{tu32} we get $E_{r1,B} , E_{r2,B} \in {\goth E}_{n,B}({\Bbb D})$. Now part (ii) of
Theorem \ref{tu51} provides the existence of a unitary matrix satisfying
$\Theta_{r2}=\Theta_{r1} \cdot V$ and $E_{r2,B}=V^{\star} \cdot E_{r1,B}$.
Hence,
\[ E_{r2}= \frac{1}{d} \cdot E_{r2,B} = \frac{1}{d} \cdot V^{\star} \cdot E_{r1,B} = V^{\star} \cdot E_{r1} . \]
Thus, (ii) is proved.\\
 Assertions (iii) and (iv) can be established analogously. \hfill
 $\boxx$\\[0.1cm]

\begin{cor}
{\sl
\label{cu51}
Let $F \in {\goth N}_n^+({\Bbb D})$ be such that $\det{F} \not\equiv 0$.
Then there exist functions $B_1 \in {\goth I}_{n,B,l}({\Bbb D}), S_1 \in
 {\goth I}_{n,s}({\Bbb D})$
and $E_1 \in {\goth E}_n({\Bbb D})$ {\rm (}resp. $B_2 \in
 {\goth I}_{n,B,r}({\Bbb D}), S_2 \in {\goth I}_{n,s}({\Bbb D})$
and $E_2 \in {\goth E}_n({\Bbb D})${\rm)} such that
\[ F = B_1 \cdot S_1 \cdot E_1 \qquad (\mbox{resp.}
 \ F = E_2 \cdot S_2 \cdot B_2)  . \]
}
\end{cor}
\ \\[-0.8cm]

 PROOF.  The assertion follows immediately by combining Theorem \ref{tu41}
and\\ Theorem \ref{tu53}. \qquad \hfill $\boxx$

It should be mentioned that using deep results and methods of V. Potapov [Pot]
an alternate approach to Theorem \ref{tu53} and Corollary \ref{cu51} was
presented by J.P. Ginzburg [Gi1]. His result contains also a multiplicative
integral representation for the outer factor and the singular inner
component.

The following theorem provides a useful characterization of the case that the
inner component in the inner - outer factorization of a given function from
${\goth N}_n^+({\Bbb D})$ is a Blaschke-Potapov product.

\begin{theo}
\label{tu54}
{\sl
Let $F \in {\goth N}_n^+({\Bbb D})$ be such that $\det{F} \not\equiv 0$.
Suppose that the functions $\Theta_r, \Theta_l \in {\goth I}_n({\Bbb D})$ and
$E_r, E_l \in {\goth E}_n({\Bbb D})$ satisfy $\Theta_r \cdot E_r = E_l \cdot \Theta_l = F$.\\
Then the following statements are equivalent:\\[0.2cm]
\hspace*{0.95cm}{\rm (i) } \parbox[t]{14.0cm}{$\quad \ \Theta_r \in {\goth I}_{n,B,r}({\Bbb D})$ }\\[0.2cm]
\hspace*{0.90cm}{\rm (ii) } \parbox[t]{14.0cm}{$\quad \ \Theta_l \in {\goth I}_{n,B,l}({\Bbb D})$}\\[0.2cm]
\hspace*{0.8cm}{\rm (iii) } \parbox[t]{14.0cm}{ \quad $\lim\limits_{s \to 1-0} \int\limits_{{\Bbb T}} \ln{|\det{[F(st)]}|} \ m(dt) = \int\limits_{{\Bbb T}} \ln{| \underline{F}(t) |} \ m(dt) $.}
}\\[0.2cm]
\end{theo}

 PROOF.  In view of the fact that $E_r, E_l \in {\goth E}_n({\Bbb D})$,
 the functions $\det{E_r}$ and $\det{E_l}$ are outer.
 Moreover, since $\Theta_r, \Theta_l \in {\goth I}_n({\Bbb D})$,
part (b) of Theorem \ref{tu42} implies that the functions $\det{\Theta_r}, \det{\Theta_l}$
are inner. From part (d) of Theorem \ref{tu42} it follows that (i) (resp. (ii))
holds  if and only if $\det{\Theta_r}$ (resp. $\det{\Theta_l}$) is a
Blaschke product. According to Lemma \ref{lu44} this is equivalent to
\begin{equation} \lim\limits_{s \to 1-0} \int\limits_{{\Bbb T}} \ln{|\det{[\Theta_r(st)]}|} \ m(dt) = 0 \label{5.12}
\end{equation}
(resp.
\begin{equation} \lim\limits_{s \to 1-0} \int\limits_{{\Bbb T}} \ln{|\det{[\Theta_l(st)]}|} \ m(dt) = 0 ). \label{5.13}
\end{equation}
From the multiplicative decomposition $F=\Theta_r \cdot E_r$ (resp. $F=E_l \cdot \Theta_l$)
it follows immediately that (\ref{5.12}) (resp. (\ref{5.13})) is equivalent
to (iii).\\ Thus, the statements (i) - (iii) are equivalent.  \hfill $\boxx$

\begin{rem}
\label{ru55}
It is instructive to compare statement (iii) in \mbox{Theorem \ref{tu54}}
with the inequality (\ref{1.9}) which is fulfilled for an arbitrary function $F$
from ${\goth N}_n^+({\Bbb D})$.
\end{rem}\\[0.3cm]
%%%%%%%%%%%%%%%%%%%%%%%%%%%%%%%%%%%%%%%%%%%%%%%%%%%%%%%%%%%%%%%%%%%%%%%%%%%%%%
%%%%%%%%%%%%%%%%%%%%%%%%%%%%%%% SECTION 6 %%%%%%%%%%%%%%%%%%%%%%%%%%%%%%%%%%%%
%%%%%%%%%%%%%%%%%%%%%%%%%%%%%%%%%%%%%%%%%%%%%%%%%%%%%%%%%%%%%%%%%%%%%%%%%%%%%%

\begin{minipage}{15.0cm}
\section{\hspace{-0.4cm}. AN ANALOGUE OF FROSTMAN'S THEOREM
 FOR  MATRIX FUNCTIONS OF THE SMIRNOV CLASS}
\end{minipage}\\[0.1cm]
\setcounter{equation}{0}

Let $f$ be a nonconstant function from the Smirnov class ${\goth N}^+({\Bbb D})$.
For $\lambda \in {\Bbb C}$ the function
\begin{equation} f_{\lambda} := f - \lambda \label{6.1}
\end{equation}
clearly belongs to ${\goth N}^+({\Bbb D})$ too. Thus, there exists an inner
 function
$\theta_{\lambda}$ and an outer function $e_{\lambda}$ such that
\begin{equation} f_{\lambda} = \theta_{\lambda} \cdot e_{\lambda} . \label{6.2}
\end{equation}
It will turn out that in some sense "the typical situation" corresponds
to the case that the inner function $\theta_{\lambda}$ in (\ref{6.2}) is a
Blaschke product. The set of all $\lambda \in {\Bbb C}$ for which
 $\theta_{\lambda}$
is not a Blaschke product is very thin. (A remarkable result of this type goes
 back to Frostman [Fr].) The corresponding notion of thinness can be
 formulated in terms of potential theory. For this reason, now we recall some
 notions of potential theory.

Suppose that $\nu$ is a nonnegative Borel measure with compact support. For
all $\xi \in {\Bbb C}$ the integral
\begin{equation} U^{(\nu)}(\xi)
 := \int\limits_{{\Bbb C}} \ln{|\xi - \lambda|} \ \nu(d\lambda) \label{6.3}
\end{equation}
is then well-defined and takes its values in $[-\infty,\infty)$.
The function $U^{(\nu)} : {\Bbb C} \to [-\infty,\infty)$ is called
the {\sl logarithmic potential of $\nu$.}
A Borel measure $\nu$ on ${\Bbb C}$ is said to be nontrivial if it is not the
zero measure. If $K$ is a Borel subset of ${\Bbb C}$, the Borel measure $\nu$
is said to be concentrated on $K$ if $\nu({\Bbb C} \setminus K) = 0$.
By definition, a Borel subset $K$ on ${\Bbb C}$ is called thin if for each
nontrivial Borel measure $\nu$ which is concentrated on $K$ the associated
logarithmic potential $U^{(\nu)}$ is not bounded from below, or in other words
if
\[ \inf\limits_{\xi \in {\Bbb C}} U^{(\nu)}(\xi) = - \infty . \]
If $K$ is not thin, then there exists a nontrivial Borel measure $\nu$ which is
concentrated on $K$ and satisfies
\begin{equation} \inf\limits_{\xi \in {\Bbb C}} U^{(\nu)}(\xi) > - \infty .
\label{6.4}
\end{equation}
The notion of {\sl logarithmic capacity} is introduced in potential theory.
More precisely, this means that with each Borel subset $K$ of ${\Bbb C}$
there is associated  a nonnegative
number $\mbox{cap} K$ which is called \mbox{{\sl the logarithmic capacity of
 $K$.}}  It turns out that a Borel subset $K$ of ${\Bbb C}$ is thin if and only
if $\mbox{cap} K = 0$. In other words, if $\mbox{cap} K > 0$, then there
exists a nontrivial Borel measure $\nu$
which is concentrated on $K$ and satisfies condition (\ref{6.4}).
If $\mbox{cap} K > 0$, then amongst all the nontrivial Borel measures
$\nu$ which are concentrated on $K$ and
satisfy (\ref{6.4}) there is a distinguished probability measure
 $\nu_{\raisebox{-0.1cm}{$\scriptstyle K$}}$,
the so-called {\sl equilibrium measure of $K$.} This measure
 $\nu_{\raisebox{-0.1cm}{$\scriptstyle K$}}$
 is a
solution of several natural extremal problems. (If $\mbox{cap} K = 0$ the
equilibrium measure is not defined.)\\
The logarithmic potential is not always continuous on ${\Bbb C}$ but only upper
semicontinuous on ${\Bbb C}$. More precisely, for all $\xi \in {\Bbb C}$,
\[ \mathop{\overline{\lim}}\limits_{\xi' \to \xi} U^{(\nu)}(\xi') \leq U^{(\nu)}(\xi) . \]
Although it is bounded below, the logarithmic potential of the equilibrium
measure need not  be continuous on ${\Bbb C}$. If the set $K$ is ``bad'' there
are so-called irregular points. Nevertheless it can be proved
(see de la Vall\'{e}e Poussin [LVP2], [LVP3]) that if $\mbox{cap} K > 0$, then
there exists a nontrivial nonnegative measure which is concentrated on $K$ and
for which the associated logarithmic potential is continuous on ${\Bbb C}$
(as already mentioned, the equilibrium measure
 $\nu_{\raisebox{-0.1cm}{$\scriptstyle K$}}$
 can not generally be
 used for this purpose). We will not enter
into such detailed and rather delicate potential - theoretical considerations.
To avoid them  we give the following definition.

\begin{defi}
\label{du61}
{\sl
{\rm A bounded Borel subset $K$ of ${\Bbb C}$} is  said {\rm to have positive
logarithmic capacity} if there exists a nontrivial Borel measure $\nu$ which is
concentrated on $K$ and for which the associated logarithmic potential is
 continuous on ${\Bbb C}$.
}
\end{defi}

Clearly, if $K_1 \subseteq K_2$ and $K_1$ is a set of positive logarithmic
capacity, then $K_2$ is also a set of positive logarithmic capacity.

{\rm LEMMA ON THE CAPACITY OF AN INTERVAL. } {\sl Every interval of the
 complex plane is
 a set of positive logarithmic capacity.}

 PROOF.   Without loss of generality we can assume that the considered
interval is a subinterval $(\alpha,\beta)$ of the real axis where
$-\infty < \alpha < \beta < \infty$. Now we take for $\nu$ the restriction
of one dimensional Lebesgue measure to this interval $(\alpha,\beta)$.
The function $U^{(\nu)} : {\Bbb C} \to [-\infty,\infty)$ which is defined by
the rule
\[ \xi \to \int\limits_{(\alpha,\beta)} \ln{| \xi - \lambda |} \
 \nu (d\lambda) \]
is continuous in ${\Bbb C}$. This can be checked in several ways, e.g. one can
 compute explicitly  and then  obtain the continuity of $U^{(\nu)}$ by direct
 estimates. \qquad \hfill $\boxx$

W. Rudin [Ru1] (see also section 3.6 of the monograph [Ru2]) proved the
following fact which generalizes Frostman's original result:\\
Let $f \in {\goth N}^+({\Bbb D})$ with $f \not\equiv 0$ and let $K$ be some
bounded Borel subset of ${\Bbb C}$ with positive logarithmic capacity.
Then there exist a $\lambda \in K$ such that the inner factor in the
 multiplicative decomposition (\ref{6.2}) is a Blaschke
product. (Indeed, W. Rudin obtained a more general result which is formulated
for the Smirnov class ${\goth N}^+({\Bbb D}^p)$ in the polydisc ${\Bbb D}^p$.
This class is a natural analogue of ${\goth N}^+({\Bbb D})$ and coincides with it
in the case $p=1$.) It should be mentioned that \mbox{S.A. Vinogradov} [Vin]
independently obtained such a generalization of Frostman's theorem too.

\begin{rem}
\label{ru61}
Let $F \in {\goth N}_n^+({\Bbb D})$ and  define
$F_{\lambda} := F - \lambda \cdot I_n$ for $\lambda \in {\Bbb C}$. Then the
set $M_F:=\{ \lambda \in {\Bbb C} \ : \ \det (F_{\lambda}) \equiv 0 \}$
is  finite.
\end{rem}

Now we formulate our main result.\\[0.2cm]

\begin{theo}
\label{tu61}
{\sl
Let $F \in {\goth N}_n^+({\Bbb D})$. Assume that for $\lambda \in {\Bbb C} \setminus M_F$
the functions $\Theta_{\lambda,r} \in {\goth I}_n({\Bbb D})$
 and $E_{\lambda,r} \in {\goth E}_n({\Bbb D})$
are factors in the multiplicative decomposition
\[ F_{\lambda} = \Theta_{\lambda,r} \cdot E_{\lambda,r} . \]
Suppose that $K$ is a bounded Borel subset of ${\Bbb C}$ with positive logarithmic
capacity. Then there exists a point
 $\lambda \in K \cap ({\Bbb C} \setminus M_F)$ for
which $\Theta_{\lambda,r}$ is a Blaschke-Potapov product.
}
\end{theo}

\begin{cor}
\label{cu61}
{\it
The set of all $\lambda \in {\Bbb C} \setminus M_F$ for which $\Theta_{\lambda,r}$
is a Blaschke-Potapov product is dense in ${\Bbb C}$.
}
\end{cor}

 PROOF :  Combine Theorem \ref{tu61} and the Lemma on the capacity
of an interval. \quad \hfill $\boxx$

In order to to follow the strategy of W. Rudin's proof we shell need to
introduce a number of classes of scalar functions of several variables.\\[0.1cm]

\begin{defi}
\label{du62}
{\sl
A function $\sigma : {\Bbb C}^n \to {\Bbb R}$ is called {\rm symmetric} if for all
permutations $\left( \begin{array}{ccc} 1 & \ldots & n\\ i_1 & \ldots & i_n \end{array} \right)$
and all $x=(x_1, \ldots, x_n)^{\top} \in {\Bbb C}^n$
the relation
\[ \sigma((x_{i_1}, \ldots, x_{i_n})^{\top}) = \sigma((x_1, \ldots, x_n)^{\top}) \]
is valid.
}
\end{defi}

In view of Definition \ref{du62} the following object is well-defined.\\[0.1cm]

\begin{defi}
{\sl
\label{du63}
Let $\sigma : {\Bbb C}^n \to {\Bbb R}$ be a symmetric function. Then the
function $\varphi_{\sigma} : {\goth M}_n \to {\Bbb R}$ which is defined by the
rule
\[ A \to \sigma((l_1(A), \ldots, l_n(A))^{\top}), \]
where $(l_j(A))_{j=1}^n$ are the roots of the characteristic polynomial of $A$
(taking into account their algebraic multiplicities),
is called {\it the function of matrix argument which is generated by the
symmetric function $\sigma$.}
}\\[0.1cm]
\end{defi}

\begin{lemma}
\label{lu61}
{\sl
Suppose that $\sigma : {\Bbb C}^n \to {\Bbb R}$ is a continuous symmetric function. Then $\varphi_{\sigma}$ is a continuous function.\\[-0.4cm]
}
\end{lemma}

 PROOF.  The lemma is an immediate consequence of Theorem 5.1 from Chapter II
in Kato's monograph [Ka]. (See there especially formula (5.3) and the text
following it.)  \hfill $\boxx$

If the symmetric function $\sigma : {\Bbb C}^n \to {\Bbb R}$ is a polynomial
or a rational function in $n$ variables $x_1, \ldots, x_n$, then it can be
expressed as a polynomial or a rational function of the elementary symmetric
functions. In this case the function $\varphi_{\sigma}$ is a polynomial or a
rational function of the elements of the matrix variable.

We introduce now a {\sl potential of the matrix argument}. Roughly speaking,
we insert a matrix argument in formula (\ref{6.3}) instead of the complex
variable.

\begin{rem}
\label{ru62}
Suppose that $A \in {\goth M}_n$. Then the function
$h_A: {\Bbb C} \to {\Bbb R}$
which is defined by $\lambda \to |\det{[A-\lambda I_n]}|$ is continuous.
Hence, the function $\ln{h_A}$ is continuous and locally bounded above. If
$\nu$ is a finite Borel measure on ${\Bbb C}$ with compact support, then the
function $\Phi^{(\nu)}: {\goth M}_n \to [-\infty,\infty)$ with
\begin{equation} \Phi^{(\nu)}(A) := \int\limits_{{\Bbb C}} \ln{|\det{(A-\lambda I_n)}|} \ \nu(d\lambda) \label{6.5}
\end{equation}
is well-defined.
\end{rem}

\begin{defi}
\label{du64}
{\sl
Suppose that $\nu$ is a finite Borel
 measure on ${\Bbb C}$ with compact support.
Then the function $\Phi^{(\nu)}: {\goth M}_n \to [-\infty,\infty)$ which is
defined by (\ref{6.5}) is called {\it the potential of the matrix argument
 associated with $\nu$.}
}
\end{defi}

Assume that $\nu$ is a finite Borel measure on ${\Bbb C}$ with compact
support. Let $U^{(\nu)}$ denote the logarithmic potential of $\nu$.
 Let $A \in {\goth M}_n$
and let $(l_k(A))_{k=1}^n$ be  roots of the characteristic
polynomial of $A$. For $\lambda \in {\Bbb C}$ we then have
\begin{equation} \ln{|\det{(A-\lambda I_n)}|} = \sum\limits_{k=1}^n \ln{|l_k(A)-\lambda|} . \label{6.6}
\end{equation}
Hence, upon taking (\ref{6.3}) into account we get
\begin{equation} \Phi^{(\nu)}(A) = \sum\limits_{k=1}^n U^{(\nu)}(l_k(A)) . \label{6.7}
\end{equation}
We define $\sigma^{(\nu)} : {\Bbb C}_n \to [-\infty,\infty)$ via
\begin{equation} \left( \begin{array}{c} x_1 \\ \vdots \\ x_n \end{array} \right) \to \sum\limits_{k=1}^n U^{(\nu)}(x_k) . \label{6.8}
\end{equation}
Obviously, the function $\sigma^{(\nu)}$ is symmetric. From Definition \ref{du63},
(\ref{6.7}) and (\ref{6.8}) we infer that
\begin{equation} \Phi^{(\nu)}(A) = \varphi_{\sigma^{(\nu)}}(A) . \label{6.9}
\end{equation}

\begin{lemma}
\label{lu62}
{\sl
Suppose that $\nu$ is a finite Borel measure on ${\Bbb C}$ with compact support
such that the associated logarithmic potential $U^{(\nu)}$ is continuous on
${\Bbb C}$. Then the function $\Phi^{(\nu)}: {\goth M}_n \to [-\infty,\infty)$
which is defined by (\ref{6.5}) is continuous on ${\goth M}_n$.
}
\end{lemma}

 PROOF.  Indeed, from (\ref{6.8}) it follows that $\sigma^{(\nu)}$ is a
continuous function on ${\Bbb C}^n$. Then in view of (\ref{6.9}) and Lemma \ref{lu61} the assertion follows.  \hfill $\boxx$

\begin{defi}
\label{du65}
{\sl
Let $\nu$ be a finite Borel measure on ${\Bbb C}$ with compact support.
Then the functions $\Phi_+^{(\nu)}: {\goth M}_n \to [0,\infty)$ and
$\Phi_-^{(\nu)}: {\goth M}_n \to [0,\infty)$ are defined via the formulas
\begin{equation} \Phi_+^{(\nu)}(A) := \int\limits_{{\Bbb C}} \ln^+{|\det{(A-\lambda I_n)}|} \ \nu(d\lambda) \label{6.10}
\end{equation}
and
\begin{equation} \Phi_-^{(\nu)}(A) := \int\limits_{{\Bbb C}} \ln^-{|\det{(A-\lambda I_n)}|} \ \nu(d\lambda) , \label{6.11}
\end{equation}
respectively.\\[0.0cm]
}
\end{defi}

\begin{lemma}
\label{lu63}
{\sl
Suppose that $\nu$ is a finite Borel measure on ${\Bbb C}$ with compact
support. Then the function $\Phi_+^{(\nu)}$  defined by (\ref{6.10})
is continuous on ${\goth M}_n$.
}
\end{lemma}

 PROOF.  The function $f: {\goth M}_n \times {\Bbb C} \to [0,\infty)$
which is defined by
\[ f(A,\lambda) := |\det{(A-\lambda I_n)}| \]
is continuous on ${\goth M}_n \times {\Bbb C}$. Since the function
$\ln^+ := \max\{ \ln , 0 \}$ is continuous on $[0,\infty)$ the composition mapping $\ln^+ f$
is continuous on ${\goth M}_n \times {\Bbb C}$. From this we infer that the function
$\Phi_+^{(\nu)}$ is continuous on ${\goth M}_n$. \qquad \hfill $\boxx$

\begin{lemma}
\label{lu64}
{\sl
Suppose that $\nu$ is a finite Borel measure on ${\Bbb C}$ with compact support
such that the associated logarithmic potential $U^{(\nu)}$ is continuous on ${\Bbb C}$.
Then the function $\Phi_-^{(\nu)}$ which is defined by (\ref{6.11}) is
continuous on ${\goth M}_n$; it is also bounded:
\begin{equation} \sup\limits_{A \in {\goth M}_n} \Phi_-^{(\nu)}(A) < +\infty . \label{6.12}\end{equation}
}
\end{lemma}

 PROOF.  From Definitions \ref{du64} and \ref{du65} we get the identity
\begin{equation} \Phi^{(\nu)} = \Phi_+^{(\nu)} - \Phi_-^{(\nu)} . \label{6.13}
\end{equation}
In view of Lemma \ref{lu62} the function $\Phi^{(\nu)}$ is continuous whereas
Lemma \ref{lu63} provides the continuity of $\Phi_+^{(\nu)}$. Thus, (\ref{6.13})
shows the continuity of $\Phi_-^{(\nu)}$. We define the functions
$U_+^{(\nu)} : {\Bbb C} \to [0,\infty)$ and $U_-^{(\nu)} : {\Bbb C} \to [0,\infty)$
by
\begin{equation} U_+^{(\nu)}(\xi) := \int\limits_{{\Bbb C}} \ln^+{|\xi-\lambda|} \ \nu(d\lambda) \label{6.14}
\end{equation}
and
\begin{equation} U_-^{(\nu)}(\xi) := \int\limits_{{\Bbb C}} \ln^-{|\xi-\lambda|} \ \nu(d\lambda). \label{6.15}
\end{equation}
Combining (\ref{6.3}), (\ref{6.14}) and (\ref{6.15}) we see that
\begin{equation} U^{(\nu)} = U_+^{(\nu)} - U_-^{(\nu)} . \label{6.16}
\end{equation}
Since the function $U^{(\nu)}$ is continuous by assumption and since the
 function
$U_+^{(\nu)}$ is always continuous (by Lemma \ref{lu63} with $n=1$) the
continuity of $U_-^{(\nu)}$ follows from (\ref{6.16}). If $(r_k)_{k=1}^n$ is
a sequence from $[0,\infty)$, then clearly
\begin{equation} \ln^-{\left( \prod\limits_{k=1}^n r_k \right)} \leq \sum\limits_{k=1}^n \ln^-{r_k} . \label{6.17}
\end{equation}
Let $A \in {\goth M}_n$ and let $(l_k(A))_{k=1}^n$ be the roots of the
characteristic polynomial of $A$. In view of (\ref{6.6})  we get
\begin{equation} \ln^-{|\det{(A-\lambda I_n)}|} = \ln^-{\left( \prod\limits_{k=1}^n |l_k(A) - \lambda| \right)} . \label{6.18}
\end{equation}
From (\ref{6.17}), (\ref{6.18}) and (\ref{6.15}) we infer that
\[ \Phi_-^{(\nu)}(A) \leq \sum\limits_{k=1}^n U_-^{(\nu)}(l_k(A)) . \]
Hence,
\begin{equation} \sup\limits_{A \in {\goth M}_n} \Phi_-^{(\nu)}(A) \leq n \cdot \sup\limits_{\xi \in {\Bbb C}} U_-^{(\nu)}(\xi) . \label{6.19}
\end{equation}
Now it remains to prove that our assumptions ensure that
\begin{equation} \sup\limits_{\xi \in {\Bbb C}} U_-^{(\nu)}(\xi) < \infty \label{6.20}
\end{equation}
is fulfilled. If $\xi \in {\Bbb C}$ satisfies
\begin{equation}
|\xi| \geq 1 + \sup\limits_{\lambda \in \mbox{{\scriptsize supp }} \nu}
 |\lambda|, \label{6.21}
\end{equation}
then using (\ref{6.15}) we see that
\begin{equation} U_-^{(\nu)}(\xi) = 0 . \label{6.22}
\end{equation}
Now the continuity of $U_-^{(\nu)}$, (\ref{6.21}), (\ref{6.22}) and a classical
theorem due to Weierstrass yield (\ref{6.20}). The lemma is proved.
\hfill $\boxx$

\begin{rem}
\label{ru63}
If $a,b \in [0,\infty)$, then
\[ \ln^+{(a+b)} \leq \ln^+{a} + \ln^+{b} + \ln{2} . \]
\end{rem}

\begin{rem}
\label{ru64}
Let $A \in {\goth M}_n$. Then \ $|\det{A}| \leq \| A \|^n$.
\end{rem}

\begin{rem}
\label{ru65}
Let $A \in {\goth M}_n$ and $\lambda \in {\Bbb C}$. Then
\[ \ln^+{|\det{[A-\lambda I_n]}|} \leq n \cdot [ \ln^+{\| A \|} + \ln^+{|\lambda|} + \ln{2}] . \]
\end{rem}\\[0.2cm]
\ \ \  Indeed, using remarks \ref{ru64} and \ref{6.3} we obtain
\begin{eqnarray} \ln^+{|\det{[A-\lambda I_n]}|} & \leq & \ln^+{[\| A-\lambda I_n \|^n]} = n \cdot \ln^+{[\| A-\lambda I_n \|]} \nonumber\\
& \leq & n \cdot \ln^+{[\| A \| + \| \lambda I_n \|]} = n \cdot \ln^+{[\| A \| + | \lambda |]} \nonumber\\
& \leq & n \cdot [ \ln^+{\| A \|} + \ln^+{|\lambda|} + \ln{2}] .
 \nonumber
\end{eqnarray}

 PROOF OF THEOREM \ref{tu61}.  Let $\lambda \in {\Bbb C}$. For $r \in [0,1)$
 we define
\begin{equation} v_r(\lambda) := \int\limits_{{\Bbb T}} \ln{|\det{[F(rt)-\lambda I_n]}|} \ m(dt) . \label{6.23}
\end{equation}
Assume that $r_1,r_2 \in [0,1)$ satisfy $r_1 \leq r_2$. Since the function
$\det{[F-\lambda I_n]}$ is holomorphic we get $v_{r_1}(\lambda) \leq v_{r_2}(\lambda)$.
Thus, the limit
\begin{equation} v_{1-0}(\lambda) := \lim\limits_{r \to 1-0} v_r(\lambda)
 \label{6.24}
\end{equation}
exists. Define
\begin{equation} v(\lambda) := \int\limits_{{\Bbb T}} \ln{|\det{[\underline{F}
(t)-\lambda I_n]}|} \ m(dt) . \label{6.25}
\end{equation}
If we apply inequality (\ref{1.9}) to the function $F-\lambda I_n$, then using
(\ref{6.23}) - (\ref{6.25}) we obtain
\begin{equation} v_{1-0}(\lambda) \leq v(\lambda) . \label{6.26}
\end{equation}
According to Theorem \ref{tu54}, equality holds  in (\ref{6.26}) for those
and only those $\lambda \in {\Bbb C} \setminus M_F$ for which the inner factor
$\Theta_{\lambda,r}$ is a Blaschke-Potapov product. Consequently, Theorem
\ref{tu54}
reduces the question which is discussed in Theorem \ref{tu61} to the study of
the structure of the set of all $\lambda \in {\Bbb C} \setminus M_F$ for which
the inequality in (\ref{6.26}) is strict. More formally, we will show that if
$K$ is a bounded Borel subset of positive logarithmic capacity then there
exists a point $\lambda \in K \cap ( {\Bbb C} \setminus M_F )$ such that
equality holds true in (\ref{6.26}). Furthermore, we will show that if
$K$ is such a set and if $\nu$ is a finite Borel measure on ${\Bbb C}$ which
is concentrated on $K$, i.e.,
\[ \nu( {\Bbb C} \setminus K ) = 0, \]
and if the associated logarithmic potential $U^{(\nu)}$
(see (\ref{6.3})) is continuous in ${\Bbb C}$, then the identity
\begin{equation} \int\limits_{{\Bbb C}} [ v(\lambda) - v_{1-0}(\lambda) ] \ \nu(d\lambda) = 0 \label{6.27}
\end{equation}
is valid. Clearly, from (\ref{6.26}) and (\ref{6.27}) it will follow that
$v(\lambda)=v_{1-0}(\lambda)$ for almost all $\lambda$ with respect to $\nu$.
In particular, there exists a $\lambda \in K \cap ({\Bbb C} \setminus M_F)$
for which
$v(\lambda)=v_{1-0}(\lambda)$ is satisfied. Now we are going to prove (\ref{6.27}).
According to (\ref{6.23}) for $r \in [0,1)$ and $\lambda \in {\Bbb C}$
we have
\begin{equation} \int\limits_{{\Bbb T}} \ln^+{|\det{[F(rt)-\lambda I_n]}|} \ m(dt) - \int\limits_{{\Bbb T}} \ln^-{|\det{[F(rt)-\lambda I_n]}|} \ m(dt) = v_r(\lambda) . \label{6.28}
\end{equation}
In view of Remark \ref{ru65}, the inequality
\begin{equation} \ln^+{|\det{[F(rt)-\lambda I_n]}|}
\leq n \cdot [ \ln^+{\|F(rt)\|} + \ln^+{|\lambda|} + \ln{2} ] \label{6.29}
\end{equation}
holds for $r \in [0,1), \lambda \in {\Bbb C}$ and $t \in {\Bbb T}$.
For $\lambda \in {\Bbb C}$ and $r \in [0,1)$ the function $G_{\lambda,r}:
 {\Bbb T} \to [0,\infty)$
is defined by
\begin{equation} G_{\lambda,r}(t) := \det{[F(rt)-\lambda I_n]} . \label{6.30}
\end{equation}
Suppose that $\lambda \in {\Bbb C}$ is fixed. Then from (\ref{6.29}) and (\ref{6.30})
we infer that the family $(\ln^+{|G_{\lambda,r}|})_{r \in [0,1)}$ is uniformly
$m$-integrable. Clearly, for almost all $t \in {\Bbb T}$ with respect to $m$
we have
\[ \lim\limits_{r \to 1-0} \ln^+{|\det{[F(rt)-\lambda I_n]}|} = \ln^+{|\det{[\underline{F}(t) - \lambda I_n]}|} . \]
Thus, using Vitali's convergence theorem again,  we get
\begin{equation}  \lim\limits_{r \to 1-0} \int\limits_{{\Bbb T}}
 \ln^+{|\det{[F(rt)-\lambda I_n]}|} \ m(dt) = \int\limits_{{\Bbb T}}
 \ln^+{|\det{[\underline{F}(t) - \lambda I_n]}|} \ m(dt) . \label{6.31}
\end{equation}
Taking into account (\ref{6.31}) we obtain the formula
\begin{eqnarray} &
 & \int\limits_{{\Bbb T}} \ln^+{|\det{[\underline{F}(t)-\lambda I_n]}|}
 \ m(dt) - \lim\limits_{r \to 1-0} \int\limits_{{\Bbb T}}
 \ln^-{|\det{[F(rt)-\lambda I_n]}|} \ m(dt) \nonumber\\
 & & = v_{1-0}(\lambda) \label{6.32}
\end{eqnarray}
by letting $r \rightarrow 1-0 $ in (\ref{6.28}),
where the limit of the second term on the left hand side of (\ref{6.32})
necessarily exists. From (\ref{6.25}) and (\ref{6.32}) it follows that
\begin{eqnarray} v(\lambda) - v_{1-0}(\lambda) & = & \lim\limits_{r \to 1-0}
 \int\limits_{{\Bbb T}} \ln^-{|\det{[F(rt)-\lambda I_n]}|} \ m(dt) \nonumber\\
 & & - \int\limits_{{\Bbb T}} \ln^-{|\det{[\underline{F}(t)-\lambda I_n]}|}
 \ m(dt) . \label{6.33}
\end{eqnarray}
In general, the family $(\ln^-{|G_{\lambda,r}|})_{r \in [0,1)}$ is not
 uniformly
$m$-integrable. For this reason, the right hand side in (\ref{6.33}) is not
necessarily zero. (However, according to Fatou's theorem this difference is
nonnegative.) Nevertheless, it will turn out that after applying the following
averaging procedure the right hand side of (\ref{6.33}) vanishes. Suppose that
$\nu$ is a finite nonnegative measure with compact support for which the
associated logarithmic potential $U^{(\nu)}$ is continuous. We will prove that
%\begin{eqnarray} & & \int\limits_{{\Bbb C}} \left( \lim\limits_{r \to 1-0}
% \int\limits_{{\Bbb T}} \ln^-{|\det{[F(rt)-\lambda I_n]}|} \ m(dt) \right.
% \nonumber\\
%& & \qquad \left. - \int\limits_{{\Bbb T}}
% \ln^-{|\det{[\underline{F}(t)-\lambda I_n]}|} \ m(dt) \right) \ \nu(d\lambda)
% = 0 . \label{6.34}
%\end{eqnarray}
\begin{equation}  \int\limits_{{\Bbb C}} \Bigg( \lim\limits_{r \to 1-0}
 \int\limits_{{\Bbb T}} \ln^-{|\det{[F(rt)-\lambda I_n]}|} \ m(dt)
  - \int\limits_{{\Bbb T}}
 \ln^-{|\det{[\underline{F}(t)-\lambda I_n]}|} \ m(dt) \Bigg) \ \nu(d\lambda)
= 0 . \label{6.34}
\end{equation}
Using Fubini's theorem and (\ref{6.11}) we get
\begin{eqnarray*}
\displaystyle
 \int\limits_{{\Bbb C}} \Bigg( \int\limits_{{\Bbb T}}
 \ln^-{|\det{[\underline{F}(t)-\lambda I_n]}|} \ m(dt) \Bigg)\, \nu(d\lambda)
  = \int\limits_{{\Bbb T}} \Bigg( \int\limits_{{\Bbb C}}
 \ln^-{|\det{[\underline{F}(t)-\lambda I_n]}|}
\, \nu(d\lambda) \Bigg) \ m(dt) \nonumber
\end{eqnarray*}
\begin{equation}
= \int\limits_{{\Bbb T}} \Phi_-^{(\nu)}(\underline{F}(t)) \ m(dt) .
 \label{6.35}
\end{equation}
In view of (\ref{6.12}) it follows that
\begin{equation} \int\limits_{{\Bbb C}} \left( \int\limits_{{\Bbb T}}
 \ln^-{|\det{[\underline{F}(t)-\lambda I_n]}|} \ m(dt) \right) \ \nu(d\lambda)
  < \infty . \label{6.36}
\end{equation}
Now we integrate identity (\ref{6.33}) with respect to $\nu$ and use
 (\ref{6.36})
to rewrite the integral of the difference as the difference of integrals.
 Then we
rewrite the second term using (\ref{6.35}) and apply Fatou's theorem to the
 first
one. Finally, we use Fubini's theorem and (\ref{6.11}) to rewrite the first
 term.
This leads us to the following estimate
\begin{eqnarray}
&\ \displaystyle   \int\limits_{{\Bbb C}} [v(\lambda) - v_{1-0}(\lambda)]
 \ \nu(d\lambda) \hspace{0.0cm}
  \nonumber\\
 = &\displaystyle \int\limits_{{\Bbb C}} \Big[ \lim\limits_{r \to 1-0}
 \int\limits_{{\Bbb T}} \ln^-{|\det{[F(rt)-\lambda I_n]}|} \ m(dt)
  - \int\limits_{{\Bbb T}}
 \ln^-{|\det{[\underline{F}(t)-\lambda I_n]}|} \ m(dt) \Big]
 \ \nu(d\lambda)
 \nonumber\\
 = &\displaystyle  \int\limits_{{\Bbb C}}
 \Big[ \lim\limits_{r \to 1-0}
 \int\limits_{{\Bbb T}} \ln^-{|\det{[F(rt)-\lambda I_n]}|} \ m(dt) \Big] \
 \nu(d\lambda)\hspace{5.5cm}\nonumber\\
 \ &\hspace{4.2cm} -\displaystyle \int\limits_{{\Bbb C}}
 \Big[ \int\limits_{{\Bbb T}} \ln^-{|\det{[\underline{F}(t)-\lambda I_n]}|}
 \ m(dt) \Big] \ \nu(d\lambda)
 \nonumber\\
 = & \displaystyle \int\limits_{{\Bbb C}}
 \Big[ \lim\limits_{r \to 1-0}\displaystyle  \int\limits_{{\Bbb T}}
\ln^-{|\det{[F(rt)-\lambda I_n]}|} \ m(dt) \Big] \ \nu(d\lambda)
 -\displaystyle  \int\limits_{{\Bbb T}} \Phi_-^{(\nu)}(\underline{F}(t))
 \ m(dt) \nonumber\\
  \leq & \mathop{\underline{\lim}}\limits_{r \to 1-0}
\displaystyle  \int\limits_{{\Bbb C}} \Big[ \int\limits_{{\Bbb T}}
 \ln^-{|\det{[F(rt)-\lambda I_n]}|} \ m(dt) \Big] \ \nu(d\lambda)
 -\displaystyle  \int\limits_{{\Bbb T}}
 \Phi_-^{(\nu)}(\underline{F}(t)) \ m(dt) \nonumber\\
\displaystyle  = & \mathop{\underline{\lim}}\limits_{r \to 1-0}
\displaystyle  \int\limits_{{\Bbb T}} \Big[\displaystyle  \int\limits_{{\Bbb C}}
 \ln^-{|\det{[F(rt)-\lambda I_n]}|}
 \ \nu(d\lambda) \Big] \ m(dt)
 -\displaystyle  \int\limits_{{\Bbb T}} \Phi_-^{(\nu)}(\underline{F}(t)) \ m(dt) \nonumber\\
  = & \mathop{\underline{\lim}}\limits_{r \to 1-0}
\displaystyle  \int\limits_{{\Bbb T}} \Phi_-^{(\nu)}(F(rt)) \ m(dt) -
 \int\limits_{{\Bbb T}} \Phi_-^{(\nu)}(\underline{F}(t)) \ m(dt).
\hspace{5.2cm}  \label{6.37}
\end{eqnarray}
According to Lemma \ref{lu64} and our choice of $\nu$, the function
 $\Phi_-^{(\nu)}$
is continuous. Thus, for almost all $t \in {\Bbb T}$ with respect to $m$ we get
\begin{equation} \lim\limits_{r \to 1-0} \Phi_-^{(\nu)}(F(rt)) =
 \Phi_-^{(\nu)}(\underline{F}(t)) . \label{6.38}
\end{equation}
Since the function $\Phi_-^{(\nu)}$ is also bounded (see Lemma \ref{lu64}),
Lebesgue's theorem on dominated convergence and (\ref{6.38}) guarantee
\begin{equation} \lim\limits_{r \to 1-0} \int\limits_{{\Bbb T}}
 \Phi_-^{(\nu)}(F(rt)) \ m(dt) = \int\limits_{{\Bbb T}}
 \Phi_-^{(\nu)}(\underline{F}(t)) \ m(dt) . \label{6.39}
\end{equation}
Thus, combining (\ref{6.37}) and (\ref{6.39}) we obtain (\ref{6.27}).\\
As explained above this completes the proof.  \hfill $\boxx$

 COMMENTS ON THEOREM \ref{tu61}.  These comments are intended to clarify
 the function-theoretic content of Theorem
\ref{tu61}. Let $t \in {\Bbb T}$. Then the
function $G: {\Bbb C} \to [-\infty,\infty)$ defined by
\[ G(\lambda) := \ln{|\det{[\underline{F}(t) - \lambda I_n]}|} \]
is subharmonic. Let $r \in [0,1)$ and let the function $G_r : {\Bbb C} \to [-\infty,\infty)$
be defined by $G_r(\lambda) := G_{\lambda,r}(t)$, where $G_{\lambda,r}$ is
given in (\ref{6.30}). Then $G_r$ is subharmonic too. From standard theorems
on integrating parametric families of subharmonic functions (see. e.g.
Ronkin [Ron, Ch.I, \S5] or Lelong and Gruman [LG, Appendix I, Proposition
 I.14])
it follows that the function $v$ defined in (\ref{6.25}) is subharmonic
and that for each $r \in [0,1)$ the function $v_r$ defined in (\ref{6.23})
is also subharmonic. Since the family $(v_r)_{r \in [0,1)}$ increases
monotonically with $r$, the function $v_{1-0}$ defined in (\ref{6.24}) is the
 upper
envelope of this family. The function $v_{1-0}$ is not necessarily subharmonic
but its regularization $v_{1-0}^* : {\Bbb C} \to [-\infty,\infty)$ which
is defined by
\begin{equation}
 v_{1-0}^*(\lambda) := \mathop{\overline{\lim}}\limits_{\lambda' \to \lambda}
 v_{1-0}(\lambda') \label{6.40}
\end{equation}
turns out to be subharmonic (see Ronkin [Ron, Ch.I, \S5],
Lelong and Gruman [LG, Appendix I, Proposition I.25], Cartan [Car]).
 Clearly, for $\lambda \in {\Bbb C}$
the inequality
\begin{equation} v_{1-0}(\lambda) \leq v_{1-0}^*(\lambda) \label{6.41}
\end{equation}
holds. According to an ingenious theorem of H. Cartan (see e.g. Ronkin
[Ron, Ch.I, \S5] or Cartan [Car]), the upper envelope of a family of
 subharmonic functions
coincides with its regularization everywhere  except for  a set of
 logarithmic
capacity zero. (For the exact formulation of Cartan's theorem and
a  proof we refer to Ronkin [Ron, Ch.I, \S5]). In particular,
\begin{equation} \mbox{cap } ( \{ \lambda : v_{1-0}^*(\lambda) >
 v_{1-0}(\lambda) \} ) = 0 . \label{6.42}
\end{equation}
Since the function $v$ is upper semicontinuous (i.e., for $\lambda \in
 {\Bbb C}$,
the inequality $v(\lambda) \geq \mathop{\overline{\lim}}\limits_{\lambda'
 \to \lambda} v(\lambda')$
holds) the inequalities
\begin{equation} v_{1-0}(\lambda) \leq v_{1-0}^*(\lambda) \leq v(\lambda) .
 \label{6.43}
\end{equation}
follows for $\lambda \in {\Bbb C}$ from (\ref{6.26}).
If we establish in some way that for all $\lambda$ belonging to some dense
subset of ${\Bbb C}$ the equality $v_{1-0}(\lambda) = v(\lambda)$ holds true
then in view of (\ref{6.43}) we obtain
\begin{equation} v_{1-0}^* \equiv v . \label{6.44}
\end{equation}
The identity $v_{1-0}(\lambda) = v(\lambda)$ for all $\lambda$ belonging to
some dense subset of ${\Bbb C}$ clearly follows  from the identity
\[ \int\limits_I [ v(\lambda) - v_{1-0}(\lambda) ] \ \mu (d\lambda) = 0, \]
where $I$ is an arbitrary one dimensional interval of ${\Bbb C}$ and $\mu$
is one dimensional Lebesgue measure. The use of H. Cartan's theorem on
upper envelopes of families of subharmonic functions for proving the smallness
 (in the
sense of capacity) of exceptional sets has many traditions in the theory of
functions of one or several complex variables. The application of the
recently created complex potential theory, in particular the analogue of
H. Cartan's theorem for the upper envelope of a family of plurisubharmonic
 functions
(see Bedford and Taylor [BT , Section 7] and Sadullaev's survey paper [Sad])
enables one to derive results on families of matrix-valued functions of a
more general type,
namely on families which  depend holomorphically on $p$ variables where
$p \in {\Bbb N}$.

Finally, we turn our attention to the left version of our main result.

\begin{theo}
\label{tu61a}
{\sl
Let $F \in {\goth N}_n^+({\Bbb D})$. Assume that for $\lambda \in {\Bbb C} \setminus M_F$
the functions\\ $\Theta_{\lambda,l} \in {\goth I}_n({\Bbb D})$ and $E_{\lambda,l} \in {\goth E}_n({\Bbb D})$
are  factors in the multiplicative decomposition
\[ F_{\lambda} = E_{\lambda,l} \cdot \Theta_{\lambda,l} . \]
Suppose that $K$ is a bounded Borel
 subset of ${\Bbb C}$ with positive logarithmic
capacity. Then there exists a $\lambda \in K \cap ({\Bbb C} \setminus M_F)$ for
which $\Theta_{\lambda,l}$ is a Blaschke-Potapov product.}
\end{theo}

 PROOF.  Use Theorem \ref{tu61}, Remark \ref{ru31a} and Remark \ref{ru41a}.
 \qquad \hfill $\boxx$

\begin{cor}
\label{cu61a}
{\sl
The set of all $\lambda \in {\Bbb C} \setminus M_F$ for which $\Theta_{\lambda,l}$
is a Blaschke-Potapov product is dense in ${\Bbb C}$.
}
\end{cor}

For further matricial generalizations of the classical theorems of Frostman [Fr],
Heins [Hei] and Rudin [Ru1] we refer the reader to the papers [Gi6] and [GiTa1] - [GiTa3].\\[0.2cm]

\vspace{0.6cm}
\begin{minipage}{7.2cm}
Victor Katsnelson \\
vDepartment of Theoretical Mathematics\\
The Weizmann Institute of Science\\
Rehovot, 76100\\ Israel\\
e-mail: katze@wisdom.weizmann.ac.il
\end{minipage}
\hfill
\begin{minipage}{7.5cm}
Bernd Kirstein \\
Mathematisches Institut\\
Universit\"at Leipzig\\
Leipzig, 04109\\
Bundesrepublik Deutschland\\
e-mail: heide@mathematik.uni-leipzig.de
\end{minipage}\\[0.8cm]
{\bf AMS subject classification:} 30D50, 46E10.
\end{document}